\documentclass[12pt]{amsart}
\newcommand{\filename}{reduction.tex}   

\usepackage{latexsym, amssymb, amscd, amsfonts}

\usepackage{microtype}
\usepackage{ifthen}
\usepackage{longtable}

\newboolean{showpicture}
\setboolean{showpicture}{false}

\newboolean{bourbaki}
\setboolean{bourbaki}{true}

\newboolean{draft}
\setboolean{draft}{false}

\newboolean{shproofs}
\setboolean{shproofs}{false}

\newcommand{\np}{\medskip\noindent}

\newcommand{\point}{\vspace{3mm}\par\refstepcounter{subsection}\noindent{\bf \thesubsection.} }
\newcommand{\tpoint}[1]{\vspace{3mm}\par\refstepcounter{subsubsection}\noindent{\em #1 {\rm(}{\em \thesubsubsection}{\rm)} ---} }

\newcommand{\bpoint}[1]{\vspace{3mm}\par\refstepcounter{subsection}\noindent{\bf \thesubsection.} {\bf #1.} }

\renewenvironment{equation}{\medskip\noindent\refstepcounter{subsubsection}\makebox[0pt][l]{({\bf\thesubsubsection})}\begin{minipage}[b]{\textwidth}$$}{$$\end{minipage}\medskip\noindent}

\ifthenelse{\boolean{draft}}
{
}
{
}

\newcommand{\bpf}{\noindent {\em Proof.  }}
\newcommand{\epf}{\qed \vspace{+10pt}}

\newcommand{\Gl}{\operatorname{GL}}
\newcommand{\GL}{\Gl}

\newcommand{\Gr}{\operatorname{Gr}}                    

\renewcommand{\geq}{\geqslant}
\renewcommand{\leq}{\leqslant}

\newcommand{\st}{\,\,|\,\,}                            
\newcommand{\Osh}{{\mathcal O}}                        
\newcommand{\Esh}{{\mathcal E}}                        
\newcommand{\Ish}{{\mathcal I}}                        

\newcommand{\vv}{\underline{v}}                        
\newcommand{\vseq}{{\underline{\mathbf{v}}}}           
\newcommand{\smu}{{\!\mu}}                             

\renewcommand{\H}{\mathrm{H}}                          
\newcommand{\Ho}{\mathrm{H}_{\circ}}                   
\newcommand{\A}{\mathrm{A}}                            
\newcommand{\B}{\mathrm{B}}                            
\newcommand{\C}{\mathrm{C}}                            
\newcommand{\D}{\mathrm{D}}                            
\newcommand{\E}{\mathrm{E}}                            
\newcommand{\F}{\mathrm{F}}                            
\newcommand{\G}{\mathrm{G}}                            
\renewcommand{\P}{\mathrm{P}}                          
\newcommand{\X}{\mathrm{X}}                            
\newcommand{\N}{\mathrm{N}}                            
\newcommand{\V}{\mathrm{V}}                            
\newcommand{\U}{\mathrm{U}}                            
\newcommand{\Q}{\mathrm{Q}}                            
\newcommand{\LGn}{\mathrm{LG}_{n}}                     
\newcommand{\bt}{\,\mbox{\tiny$\boxtimes$}\,}          
\newcommand{\W}{\mathcal{W}}                           
\newcommand{\delneg}{\Delta^{\!-}}                     
\newcommand{\delpos}{\Delta^{\!+}}                     
\newcommand{\invset}{\Phi}                             
\newcommand{\wo}{w_{0}}                                
\newcommand{\fo}{f_{\circ}}                            
\newcommand{\fopush}{f_{\circ*}}                       
\newcommand{\so}{s_{\circ}}                            
\renewcommand{\c}{\operatorname{c}}                    
\newcommand{\T}{\mathrm{T}}                            
\renewcommand{\L}{\mathrm{L}}                          
\newcommand{\LC}{\CC}                                  
\renewcommand{\SS}{\mathrm{S}}                         
\newcommand{\spec}{\operatorname{Spec}}                

\newcommand{\Gtil}{\overline{\G}}                      
\newcommand{\Btil}{\overline{\B}}                      
\newcommand{\Ttil}{\overline{\T}}                      
\newcommand{\Xtil}{\overline{\X}}                      
\newcommand{\mut}{\overline{\mu}}                      
\newcommand{\smut}{\!\mut}                             
\newcommand{\mup}{\mu'}                                
\newcommand{\nut}{\overline{\nu}}                      
\newcommand{\snut}{\!\nut}                             
\newcommand{\nup}{\nu'}                                
\newcommand{\pitil}{{\pi}}                             
\newcommand{\etil}{\overline{e}}                       
\newcommand{\Ltil}{\overline{\L}}                      
\newcommand{\PI}{\P_{\I}}                              
\newcommand{\GI}{\G_{\I}}                              
\newcommand{\XI}{\X_{\I}}                              
\newcommand{\BI}{\B_{\I}}                              
\newcommand{\op}{\operatorname{op}}                    

\newcommand{\ggI}{\mathfrak{g}_{\I}}                   
\newcommand{\gbI}{\mathfrak{b}_{\I}}                   
\newcommand{\gbIp}{\mathfrak{b}_{\I}^{+}}              
\newcommand{\gttil}{\overline{\mathfrak{t}}}           
\newcommand{\gbtil}{\overline{\mathfrak{b}}}           
\newcommand{\ggtil}{\overline{\mathfrak{g}}}           
\newcommand{\gbtilp}{\overline{\mathfrak{b}}\vphantom{\mathfrak{b}}^{+}}
\newcommand{\ptil}{\overline{p}}                       

\newcommand{\ga}{{\mathfrak{a}}}
\newcommand{\gb}{{\mathfrak{b}}}

\newcommand{\gt}{\mathfrak{t}}

\newcommand{\LRC}{{\mathcal{C}}}                       
\newcommand{\die}{\partial}                            
\newcommand{\mult}{\operatorname{mult}}                

\newcommand{\Lie}{\operatorname{Lie}}                  
\renewcommand{\ggg}{\mathfrak{g}}                      

\newcommand{\ep}{\varepsilon}                          
\newcommand{\ept}{\overline{\ep}}
\newcommand{\epp}{{\ep}'}
\newcommand{\pu}{\overline{p}}

\newcommand{\quot}{/\!\!/}

\newcommand{\Span}{\operatorname{span}}                

\newcommand{\I}{\mathrm{I}}

\newcommand{\ott}{\operatornamewithlimits{\otimes}}

\ifthenelse{\boolean{bourbaki}}
{
\renewcommand{\AA}{\mathbf{A}} 
\newcommand{\CC}{\mathbf{C}} 
\newcommand{\PP}{\mathbf{P}} 
\newcommand{\QQ}{\mathbf{Q}} 
\newcommand{\ZZ}{\mathbf{Z}} 
}
{
\renewcommand{\AA}{\mathbb{A}} 
\newcommand{\CC}{\mathbb{C}} 
\newcommand{\PP}{\mathbb{P}} 
\newcommand{\QQ}{\mathbb{Q}} 
\newcommand{\ZZ}{\mathbb{Z}} 
}

\ifthenelse{\boolean{draft}}
{

\newcommand{\remind}[1]{{\bf[#1]}}
\newcommand{\lremind}[1]{{\bf[label:  #1]}}
\newcommand{\bremind}[1]{{\bf[label:  #1]}}

\newcommand{\comment}[1]{{\bf [#1]}}
}
{

\newcommand{\remind}[1]{{}}
\newcommand{\lremind}[1]{{}}
\newcommand{\bremind}[1]{{}}

\newcommand{\comment}[1]{{}}
}

\newcommand{\hiddenproof}[1]{
\ifthenelse{\boolean{shproofs}}{
\medskip
\begin{centering}
\begin{minipage}{0.9\textwidth}
\hrule
\vspace{-0.75\baselineskip}
\small
#1
\vspace{0.25\baselineskip}
\hrule
\end{minipage}\\
\end{centering}
\medskip
}
{
}
}

\begin{document}
\pagestyle{plain} \title{{ \large{Reduction rules for Littlewood-Richardson coefficients}}}
\author{Mike Roth}
\thanks{Research partially supported by an NSERC grant}


\begin{abstract}
Let $\G$ be a semisimple algebraic group over an algebraically-closed field of characteristic zero.
In this note we show that every regular 
face of the Littlewood-Richardson cone of $\G$ gives
rise to a {\em reduction rule}: a rule which, given a problem ``on that face'' 
of computing the multiplicity of an irreducible
component in a tensor product, reduces it to a similar problem on a group $\Gtil$ of smaller rank.

\np
In the type $\A$ case this result has already been proved by 
Derksen and Weyman using quivers, and by King, Tollu, and Toumazet using puzzles.
The proof here is geometric and type-independent.

\np
\end{abstract}

\maketitle
\tableofcontents

\date{\today.\hspace{0.5cm}  {\em \filename}}


\section{Introduction}

\np
This note is concerned with {\em reduction rules} --- rules reducing the problem of computing 
the multiplicity of an irreducible component in a tensor product of $\G$-representations 
to a similar problem on a group $\Gtil$ of smaller rank.  
The main result is that every regular codimension-$r$ face of the Littlewood-Richardson
cone of $\G$ gives rise to a rule reducing every problem on that face to 
a group whose rank is $r$ less than the rank of $\G$.

\np
Let $\G$ be a semisimple algebraic group over an algebraically closed field of characteristic zero.  
For a dominant weight $\mu$ and a representation $\V$ of $\G$ we denote by $\mult_{\G}(\V_{\mu},\V)$ the multiplicity
of the irreducible $\G$-representation $\V_{\smu}$ in $\V$.  

\np
For any $k\geq 2$ the {\em Littlewood-Richardson cone} $\LRC(k)$ is defined as the rational cone generated by 
$(\mu_1,\ldots, \mu_k,\mu)$ such that $\V_{\smu}$ is a component of $\V_{\smu_1}\otimes\cdots\otimes\V_{\smu_k}$.
It is known that $\LRC(k)$ is polyhedral, and minimal equations equations for $\LRC(k)$ are known through the
work of Belkale-Kumar \cite{BK1} and Ressayre \cite{Re}.  
A face of $\LRC(k)$ is called {\em regular} if it intersects the locus of strictly dominant weights.

\np
By the results in \cite{Re}, the regular faces of $\LRC(k)$ 
are described by the data of a subset $\I$ of the simple roots and elements $w_1$,\ldots, $w_k$, and $w$
of the Weyl group of $\G$ satisfying some conditions relative to $\I$ (see \eqref{eqn:Iconditions} for the exact
conditions). A point $(\mu_1,\ldots, \mu_k,\mu)\in\LRC(k)$ is on the face described by this data if and only if 
the weight $\sum_{i=1}^{k} w_i^{-1}\mu_i - w^{-1}\mu$ can be written as a 
$\QQ$-linear combination of elements in $\I$.

\np
Suppose that this last condition holds.  
Then let $\Gtil$ be the semisimple part of the parabolic subgroup $\PI$ determined by $\I$ and $\mut_1$,\ldots, $\mut_k$, and
$\mut$ be the restriction of the weights $w_1^{-1}\mu_1$, \ldots, $w_k^{-1}\mu_k$ and $w^{-1}\mu$ respectively 
to $\Gtil$ (see \S\ref{sec:GtilandGI} and the examples in \S\ref{sec:examples} 
for a more precise description of this process).  The main result of this paper is the construction
of a geometric map 
$(\Gtil/\Btil)^{k+1}\longrightarrow (\G/\B)^{k+1}$ 
such that pullback of global sections of a particular line bundle induces an isomorphism of vector spaces

$$(\V_{\smu_1}\otimes\cdots\otimes\V_{\smu_k}\otimes\V_{\smu}^{*})^{\G} \stackrel{\sim}{\longrightarrow}
(\V_{\smut_1}\otimes\cdots\otimes\V_{\smut_k}\otimes\V_{\smut}^{*})^{\Gtil}.$$

\np
Taking dimensions then gives the equality

$$\mult_{\G}(\V_{\smu},\V_{\smu_1}\otimes\cdots\otimes\V_{\smu_k}) = 
\mult_{\Gtil}(\V_{\smut},\V_{\smut_1}\otimes\cdots\otimes\V_{\smut_k}),$$

\np
yielding a reduction rule.

\np
One might guess that a reduction rule occurs because the individual weights $\mu_1$,\ldots, $\mu_k$, 
and $\mu$ are somehow themselves ``special'', e.g., somehow come from a group of smaller rank.  
However, since the faces in question are 
regular, at a general point of each face all the weights are strictly dominant, and so in some sense generic.  
It is instead the special configuration of the multiplicity problem ---
as witnessed by the location of the point $(\mu_1,\ldots, \mu_k,\mu)$ 
on the boundary of $\LRC(k)$ --- that allows the reduction.

\np
Given the data of $\I$ and $w_1$,\ldots, $w_k$, and $w$, it is easy to write out explicitly what
the corresponding reduction rule does, and examples are given in \S\ref{sec:examples}. 

\np
An elementary way to describe 
$\Gtil$ is to note that its Dynkin diagram is the full subdiagram of the Dynkin diagram for $\G$ corresponding to
the simple roots in $\I$. If the resulting subdiagram is disconnected then $\Gtil$ is a product of simple groups and
hence the reduction rule can also be interpreted as a factorization rule.
Under this name,
the main result of this note was already known in the type $\A$ case and was proved independently  
by Derksen and Weyman \cite[Theorem 7.14]{dw} using quivers and by 
King, Tollu, and Toumazet \cite[Theorem 1.4]{ktt} using
puzzles.  The proof here is geometric and type-independent.

\np
In type $\A$ the Littlewood-Richardson coefficients are also the structure constants in the cohomology
rings of the Grassmanians $\G/\P$ for maximal parabolic subgroups $\P$, and one might hope to generalize
the reduction rules for type $\A$ in this direction instead.  For results along this line, see the
forthcoming paper \cite{kp} of Kevin Purbhoo and Allen Knutson.

\np
{\bf Acknowledgements.} The idea that such reduction rules should hold occurred in joint work (\cite{dr1} and \cite{dr2})
with my colleague Ivan Dimitrov, and several of the ideas used in the proof of the main theorem were developed in 
\cite{dr1}.  I am also greatful to Ivan for valuable discussions on some aspects of the present paper.  I thank
Kevin Purbhoo for telling me about the paper of Derksen and Weyman and his work with Allen Knutson, as well as
for advice on the examples.  
Explicit instances of the examples were computed with the help of the computer program LiE \cite{lie}.

\section{Preliminary material}
\label{sec:prelim}
\bpoint{Notation and conventions}
Throughout this note we fix a semisimple connected algebraic group $\G$, a Borel subgroup $\B\subset\G$ and a
maximal torus $\T\subset\B$.  Related groups, whose definition depends on the choice of a subset $\I$ of simple
roots, are discussed in \S\ref{sec:GtilandGI}.  The Lie algebras of algebraic groups are denoted by 
fraktur letters, e.g. $\ggg$, $\gb$, $\gt$, etc.  We use the term ``weight'' both for characters of $\T$ and weights of $\gt$.  For a dominant
weight $\mu$ we denote by $\V_{\smu}$ the irreducible $\G$-representation of highest weight $\mu$.

\np
Let $\Delta$ denote the set of roots of $\G$ (with respect to $\T$).  
For any subset $\Phi\subset\Delta$ we denote by $\Span_{\ZZ}\Phi$ the set of integer combinations
of elements of $\Phi$.  Similarly, $\Span_{\QQ_{\geq 0}}\Phi$ and $\Span_{\ZZ_{\leq 0}}\Phi$ denote respectively 
the set of non-negative rational combinations and non-positive integer combinations of elements of $\Phi$.

\np
We denote the Weyl group of $\G$ by $\W$ and use $\ell(w)$ for the length of any $w\in\W$.
We are working over an algebraically closed field of characteristic zero; for notational convenience we will assume 
that the field is $\CC$.

\bpoint{Inversion sets}
Let $\delpos$ be the set of positive roots of $\ggg$ (with respect to $\B$). Following Kostant
\cite[Definition 5.10]{K1}, for any element $w$ of the Weyl group $\W$ 
we define $\invset_{w}$, the {\em inversion set} of $w$, to be the set of positive roots  sent to 
negative roots by $w$, i.e., 

$$
\invset_{w} := w^{-1}\delneg\cap \delpos.
$$

\np
For a subset $\Phi$ of $\delpos$, we set $\Phi^{\c}:=\delpos\setminus\Phi$. 
From the definition it follows easily that 
$\invset_{\wo w} = \invset_{w}^{\c}$ and that 
$w^{-1}\Delta^{+} = \invset_{w}^{\c} \sqcup -\invset_{w}$, and we will use these formulas without comment
in the rest of the note.

\bpoint{Discussion of $\GI$ and $\Gtil$}
\label{sec:GtilandGI}
Given a subset $\I$ of simple roots, let $\PI$ be the corresponding parabolic subgroup, $\GI$ 
the reductive part (i.e., the Levi component) of $\PI$, and $\Gtil$ the semi-simple part of $\PI$.
We define $\Delta_{\I}$ to be the roots of $\GI$.
Equivalently $\Delta_{\I}$ is the subset of $\Delta$ consisting of those roots in $\Span_{\ZZ}\I$.
We denote by $\Delta_{\I}^{+}$ the intersection $\Delta_{\I}\cap\Delta^{+}$, i.e., the positive roots of $\GI$.
Equivalently $\Delta_{\I}^{+}$ is the subset of $\Delta$ consisting of those roots in $\Span_{\ZZ_{\geq 0}}\I$.
As remarked in the introduction, the Lie algebra $\ggtil$
has an elementary description: the Dynkin diagram of $\ggtil$ is the complete subdiagram of the Dynkin 
diagram of $\ggg$ containing the nodes corresponding to the simple roots in $\I$.

\np
By definition, $\T\subseteq\GI$.
Let $\A$ be the connected component of the center of $\GI$.  Then $\A\subseteq\T$ and $\A\cap\Gtil$ 
is a finite group.  The natural map $\Gtil\times\A\longrightarrow\GI$ sending a pair of elements to their product
is a surjective map with finite kernel  
and thus induces an isomorphism at the level of Lie algebras.    We will need to use a specific fact 
about the resulting direct sum decomposition of $\ggI$ and so we describe this decomposition in more detail below.

\np
Let $\Ttil$ be the connected component of $\T\cap\Gtil$, so that $\Ttil$ is a maximal torus for $\Gtil$, and let
$\gttil=\Lie(\Ttil)$.
Since $\I$ is a set of simple roots for $\Gtil$, the restriction of the roots in $\I$ to $\gttil$ is
a basis (over $\CC$) of the dual of $\gttil$.  Hence, letting $\ga\subseteq\gt$ be the subalgebra annihilated
by the roots in $\I$ we obtain a direct sum decomposition $\gt=\gttil\oplus \ga$.

\np
By the definition of $\ga$ we have the following result which we record for later use:

\tpoint{Lemma} \label{lem:cancelling}
If $\gamma\in\Span_{\QQ}\I$, then the restriction of $\gamma$ to $\ga$ is zero.

\np
In particular, 
for any root $\alpha$ of $\ggtil$, any $x\in \ggtil^{\alpha}$ and $a\in \ga$ we have 
$[a,x] = \alpha(a)x = 0\cdot x = 0$ and hence the decomposition of $\gt$ extends to a direct sum decomposition 
$\ggI = \ggtil\oplus \ga$.  

\np
Setting $\BI:=\GI\cap\B$ and $\Btil:=\Gtil\cap\B$ then $\BI$ and $\Btil$ are Borel subgroups of $\GI$ and $\Gtil$
respectively.  The direct sum decomposition of $\ggI$ restricts to a decomposition $\gbI=\gbtil\oplus \ga$.
Note that $\BI$ and $\Btil$ have the same unipotent part (equivalently, $\gbI$ and $\gbtil$ have the same nilpotent 
part); the difference between the two groups being in their maximal tori.

\np
{\bf Restriction of weights.}
Given a weight $\mu$ 
and an element  $w\in \W$ we will use 
$\mut$ and $\mup$ for the restrictions of $w^{-1}\mu$ to $\gttil$ and $\ga$ respectively under the splitting 
$\gt=\gttil\oplus\ga$ above.  This notation omits the element $w\in\W$ used, but any time we use this notation we
will be careful to explicitly specify which element $w$ is meant for that particular restriction.

\np
In the reduction theorem it is implicit that if $\mu$ is dominant then the restriction $\mut$ 
will also be dominant.
For completeness, let us see why this is true.
Let $\kappa(\cdot,\cdot)$ be the Killing form, and
suppose that $w\in\W$ is such that $\invset_{w}\cap\Delta_{\I}^{+}=\emptyset$;  this hypothesis will hold for all
$w$ we use when reducing to $\gttil$.  
Since $w\Delta_{\I}^{+}\subseteq\Delta^{+}$, 
if $\mu$ is dominant with respect to $\gb$ 
then $\kappa(w^{-1}\mu,\alpha)=\kappa(\mu,w\alpha)\geq0$ for all $\alpha\in\Delta_{\I}^{+}$ 
and thus the restriction $\mut$ of $w^{-1}\mu$ to $\gttil$ is dominant with respect to $\gbtil$. 
Note that this argument also shows that the restriction of a strictly dominant weight is again strictly dominant, and
that the restriction of an integral weight is integral with respect to $\Ttil$.   For this reason we will also refer
to the process as ``restricting the weight to $\Ttil$''.

\np
{\bf Surjections and $\gbI$-invariants.}
We will need 
the following result giving a condition ensuring that a surjection of $\gbI$-modules induces
an isomorphism of $\gbI$-invariants.

\tpoint{Lemma} \label{lem:bI-invts}
Suppose that $0\longrightarrow\E_1\longrightarrow\E_2\longrightarrow\E_3\longrightarrow0$ is an 
exact sequence of $\gbI$-modules, and that no weight of $\E_1$ is contained in $\Delta_{\I}^{+}\cup\{0\}$.
Then the induced map $\E^{\gbI}_2\longrightarrow {\E}^{\gbI}_3$ of $\gbI$-invariants is an isomorphism.

\bpf
The first four terms of the long exact sequence arising from taking $\gbI$-invariants is
$$
0\longrightarrow\E_1^{\gbI}\longrightarrow\E_2^{\gbI}\longrightarrow\E_3^{\gbI}\longrightarrow\H^1(\gbI,\E_1).
$$
By hypothesis, the zero weight does not appear in $\E_1$, and hence $\E_1$ has no $\gt$-invariants, and so
no $\gbI$-invariants, i.e., $\E_1^{\gbI}=0$. Let $\gbIp=[\gbI,\gbI]$ be the nilpotent radical of $\gbI$. 
Since taking $\gt$ invariants is exact, the Hochschild-Serre spectral sequence for the cohomology of $\gbI$ degenerates and
we have $\H^{i}(\gbI,\E_1)=\H^{i}(\gbIp,\E_1)^{\gt}$ 
for all $i\geq 0$, and in particular for $i=1$.  
The degree one piece of the complex computing $\gbIp$-cohomology is $\C^1(\gbIp,\E_1)=(\gbIp)^{*}\otimes\E_1$.  
By hypothesis no weight of $\E_1$ lies in $\Delta_{\I}^{+}$, hence $\C^1(\gbIp,\E_1)$ has no $\gt$-invariants.  Since
the differential maps of the complex are $\gt$-equivariant this gives $\H^1(\gbIp,\E_1)^{\gt}=0$.
\epf

\bpoint{The Borel-Weil theorem}
Let $\X:=\G/\B$ and let $e\in\X$ be the image of $1_{\G}\in\G$ under the quotient map.
The restriction map sending a vector bundle $\Esh$ on $\X$ to its fibre $\E$ over $e\in\X$ induces
an equivalence of categories between the $\G$-equivariant bundles on $\X$ and representations of $\B$.  We will use 
the following special case of that equivalence in establishing the reduction rule:

\tpoint{Principle} \label{princip}
Let $\Esh$ be a $\G$-equivariant vector bundle on $\X$, and $\E$ the fibre over $e\in\X$.
then restriction of global sections to the fibre $\E$ induces an isomorphism 
$\H^0(\X,\Esh)^{\G}\stackrel{\sim}{\longrightarrow}\E^{\B}$.

\np
For any weight $\lambda$ we denote by $\L_{\lambda}$ the $\G$-equivariant line bundle on $\X$ corresponding
to the one-dimensional $\B$-representation $\CC_{-\lambda}$, i.e., the representation  where $\B$ 
acts through its quotient $\T$ with weight $-\lambda$.
The Borel-Weil theorem identifies the $\G$-representation $\H^{0}(\X,\L_{\lambda})$ for any weight $\lambda$.
The main step in the proof of the Borel-Weil theorem is the following result.

\tpoint{Lemma}\label{lem:Stab-Grep}
Suppose that $\L$ is a $\G$-equivariant line bundle on $\X$, $x\in \X$ any point, and $\B_{x}$ the stabilizer
subgroup of $\X$.  Using $\L_{x}$ for the fibre of $\L$ at $x$ and setting $\V=\H^{0}(\X,\L)$ 
then the $\B_{x}$-equivariant restriction map $\V\longrightarrow \L_{x}$ 
at $x$ identifies $\V$ as the unique irreducible representation of $\G$ (if one exists) which 
has a $\B_{x}$-equivariant surjection onto the one-dimensional $\B_{x}$-representation $\L_{x}$.  If no such
irreducible representation exists then $\V=0$.

\np
If $\lambda$ is dominant then one can show that $\H^{0}(\X,\L_{\lambda})\neq 0$.
Since the only surjective $\B$-equivariant quotient map from an irreducible representation
$\V$ onto a one-dimensional representation is projection is onto the lowest weight vector of $\V$, 
if $\L=\L_{\lambda}$ and $x=e$ then Lemma \ref{lem:Stab-Grep} yields the Borel-Weil theorem:

\tpoint{Theorem} 
\label{thm:Borel-Weil-X}
For any weight $\lambda$

$$\H^{0}(\X,\L_{\lambda}) = \left\{{
\begin{array}{cl}
\V_{\lambda}^{*} & \mbox{if $\lambda$ is a dominant weight} \\
0 & \mbox{otherwise.} \\
\end{array}}\right. 
$$

\np
Now let $\XI:=\GI/\BI^{\op}$.  Then $\XI$ is isomorphic to $\GI/\BI$ as a $\GI$-variety, but the image of $1_{\GI}$
under the quotient map $\GI\longrightarrow\XI$ has stabilizer $\B^{\op}$ instead of $\B$.  The only surjective
$\B^{\op}$-equivariant quotient map from an irreducible representation $\V$ onto a one-dimensional representation
is the projection onto the highest weight vector of $\V$.  Applying 
Lemma \ref{lem:Stab-Grep} and the splitting of $\ggI$ from \S\ref{sec:GtilandGI} then gives 
the following version of the Borel-Weil theorem for $\XI$:

\tpoint{Theorem} 
\label{thm:Borel-Weil-XI}
Suppose that $\L$ is a $\GI$ equivariant line bundle on $\XI$ with torus weight $\nu$ at the image of $1_{\GI}$ in 
$\XI$, and let $\nut$ and $\nup$ be the restrictions of $\nu$ to $\gttil$ and $\ga$ respectively under the splitting
from \S\ref{sec:GtilandGI}.  Then as a $\ggI$-module

$$\H^{0}(\XI,\L) = \left\{{
\begin{array}{cl}
\V_{\snut}\otimes\CC_{\nup} & \mbox{if $\nut$ is a dominant weight for $\gttil$} \\
0 & \mbox{otherwise.} \\
\end{array}}\right. 
$$

\np
Note that if $\L$ is the restriction to $\XI$ of a globally generated line bundle under some embedding
$\varphi\colon\XI\longrightarrow\X$, then $\H^{0}(\XI,\L)\neq 0$ and hence only the first alternative above applies.
This will be the case in the application of Theorem \ref{thm:Borel-Weil-XI} in 
Proposition \ref{prop:XI-embedding} below.

\bpoint{Schubert Varieties} \label{sec:Schub}
For any element $w\in\W$ of the Weyl group the {\em Schubert variety} $\X_{w}$ is defined by

$$\X_{w} := \overline{\B \dot{w}\B/\B}\subseteq \G/\B=\X$$

\np
where $\dot{w}$ is any lift of $w$ to $\G$.  Since everything we define using $w$ will be independent of the lift,
we will almost always omit mention of lifting and just use $w$ in place of $\dot{w}$. 
The one exception to this convention is Proposition \ref{prop:XI-embedding} below where we explicitly consider
the lift in order to show that the construction in the proposition is independent of the lifting.

\np
Recall that the classes of the Schubert cycles $\{[\X_{w}]\}_{w\in \W}$ give 
a basis for the cohomology ring $\H^{*}(\X,\ZZ)$ of $\X$.
Each $[\X_{w}]$ is a cycle of complex dimension $\ell(w)$.
The dual Schubert cycles $\{[\Omega_{w}]\}_{w\in \W}$, given by $\Omega_{w}:=\X_{\wo w}$, also form a basis. 
Each $[\Omega_{w}]$ is a cycle of complex codimension $\ell(w)$.

\np
{\bf Remark.} If $w_1$,\ldots, $w_k$, and $w\in\W$ are such that 
$\ell(w)=\sum \ell(w_i)$, then the intersection $\cap_{i=1}^{k}[\Omega_{w_i}]\cdot[\X_{w}]$ is a number.
This number is the coefficient
of $[\Omega_{w}]$ when writing the product $\cap_{i=1}^{k}[\Omega_{w_i}]$ in terms of the basis 
$\{[\Omega_{v}]\}_{v\in \W}$.

\np
To reduce notation we also use $w$ to refer to the point $w\B/\B\in \X_{w} \subseteq\X$. 
In particular for the identity element $e\in\W$,  $\X_{e}=\{e\}$.  
Note that $e\in\X$ is also the image of $1_{\G}$ under the projection from $\G$ onto $\X$.  

\np
{\bf Open affine cells of Schubert varieties.}
For any $v\in\W$ the variety $\U_{v}:=\B v\B/\B\subseteq\X_{v}$ is $\B$-stable open affine subset of $\X_{v}$ 
containing $v$ and isomorphic to affine space $\AA^{\!\ell(v)}$.  Since $\U_{v}$ is $\B$-stable 
its coordinate ring 
$\H^{0}(\U_{v},\Osh_{\U_{v}})$
decomposes into $\T$-eigenspaces.  Explicitly, 
$\U_{v}=\spec(\CC[z_{-\alpha}]_{\alpha\in\invset_{v^{-1}}})$ where each
$z_{-\alpha}$ is an independent variable on which $\T$ acts via the weight $-\alpha$. The origin of this affine
space corresponds to the point $v$.

\np
For a sequence $\vv=(v_1,\ldots, v_k)$ of elements of $\W$ we set $\U_{\vv}=\U_{v_1}\times\cdots\times\U_{v_k}$.
For any weight $\delta$ let $\H^0(\U_{\vv},\Osh_{\U_{\vv}})_{\delta}$ be the subspace of $\H^{0}(\U_{\vv},\Osh_{\U_{\vv}})$ of $\T$-eigenfunctions where $\T$ acts via $\delta$.  The above description of $\U_{v}$ immediately gives the
following easy result.

\tpoint{Lemma}\label{lem:zero-unless-positive}
For any sequence $\vv$, 
if $\delta\not\in\Span_{\ZZ_{\leq0}}\Delta^{+}$ then $\H^0(\U_{\vv},\Osh_{\U_{\vv}})_{\delta}=0$.

\np
We now come to the main constructions of this section.

\tpoint{Proposition} \label{prop:XI-embedding}
Let $v$ be an element of the Weyl group such that 
$\Delta_{\I}^{+}\subseteq\invset_{v^{-1}}$, 
$\dot{v}$ any lift of $v$ to $\G$, and $\Psi\colon\GI\longrightarrow \G$ the map 
defined by $\Psi(g)=g\dot{v}$ for
all $g\in\GI$.  
Then 

\begin{enumerate}
\item The image of $\GI$ under the composite map $\GI\stackrel{\Psi}{\longrightarrow} \G\longrightarrow \X$ 
is isomorphic to $\XI:=\GI/\BI^{\op}$ and induces a $\GI$-equivariant embedding $\psi_{v}\colon\XI\longrightarrow\X$,
independent of the lift $\dot{v}$ chosen (here $\GI$ acts on $\X$ through its inclusion $\GI\hookrightarrow\G$ as a subgroup of $\G$).
\medskip
\item The image of $\psi_{v}$ lies in $\X_{v}$. Setting $\U_{v}=\B\dot{v}\B/\B$ to be the $\B$-stable open 
affine space around $v\in\X_{v}$, then the ideal of $\XI|_{\U_{v}}$ is a direct sum of the $\T$-eigenspaces 
consisting of those functions on $\U_{v}$ with torus weight contained in 
$$\SS:=\left(\Span_{\ZZ_{\leq0}}(\Delta^{+}\setminus\Delta_{\I}^{+})\right)\setminus\{0\}.$$ 
\item 
Let $\varphi_{v}$ be 
the induced inclusion $\varphi_{v}\colon\XI\longrightarrow\X_{v}$ (i.e., $\psi_{v}$ considered as a map to $\X_{v}$).
For any dominant weight $\lambda$, the pullback map 
$\H^0(\XI,\varphi_{v}^{*}(\L_{\lambda}|_{\X_{v}}))\stackrel{\varphi_{v}^{*}}{\longleftarrow}\H^0(\X_{v},\L_{\lambda}|_{\X_{v}})$ 
is surjective, and $\H^{0}(\XI,\varphi_{v}^{*}(\L_{\lambda}|_{\X_{v}})) = \V_{\mut}\otimes \LC_{\mup}$ as a representation of
$\ggI$, where $\mut$ and
$\mup$ are the restrictions of $-v\lambda$ to $\gttil$ and $\ga$ respectively under the decomposition
$\gt=\gttil\oplus\ga$ from \S\ref{sec:GtilandGI}.
\end{enumerate}

\bpf
Two elements $g_1$, $g_2$ of $\GI$ have the same image under the composite map if and only if there is a $b\in\B$ 
such that $g_1\dot{v}=g_2\dot{v}b$, i.e., $g_2^{-1}g_1=\dot{v}b\dot{v}^{-1}$, or equivalently, 
if $g_1$ and $g_2$ are in the same coset
of the subgroup $\H:=\GI\cap\dot{v}\B \dot{v}^{-1}$.  Let $\Ho$ be the connected component of the identity of $\H$.
Since $\GI$ and $\dot{v}\B\dot{v}^{-1}$ both contain $\T$, $\Ho$
is determined by its torus weights on the tangent space at the identity. 
For every root $\alpha\in\Delta$ exactly one of $\pm\alpha$ is a root of $\dot{v}\B\dot{v}^{-1}$, and so $\Ho$
must be a Borel subgroup of $\GI$.  
This implies that $\H=\Ho$, since $\Ho$ is normal in $\H$ and since every Borel subgroup of $\GI$ is its own normalizer.
The roots of $\dot{v}\B\dot{v}^{-1}$ are $v\Delta^{+}= -\invset_{v^{-1}}\sqcup \invset_{v^{-1}}^{\c}$;
by hypothesis $\Delta_{\I}^{+}\subseteq\invset_{v^{-1}}$ and so $\Ho$ must contain $\BI^{\op}$. Thus $\Ho=\BI^{\op}$
and the image of $\GI$ under the composite map is $\XI$.
The induced map $\psi_{v}$
is independent of the lift of $v$ since $\T\subseteq\GI$, and it is clear from the description that $\psi_{v}$ is 
$\GI$-equivariant.  This proves ({\em a}). 

\np
Let $\U_{v}$ be the affine space $\B\dot{v}\B/\B$.  
Under the composite map from $\GI$ to $\X$ inducing $\psi_{v}$, the image 
$\U_{\I,v}:=\BI\dot{v}\B/\B$ of $\BI$
forms an open cell of $\psi_{v}(\XI)$ around $v\in\psi_{v}(\XI)$.  
Since $\BI\subseteq\B$ this shows that $\U_{\I,v}$ is contained in $\U_{v}$ and hence, taking Zariski closures in 
$\X$, that $\psi_{v}(\XI)$ is contained in $\X_{v}$.

\np
By the above discussion on open affine cells, 
$\U_{v}=\spec(\CC[z_{-\alpha}]_{\alpha\in\invset_{v^{-1}}})$ where each
$z_{-\alpha}$ is an independent variable on which $\T$ acts via the weight $-\alpha$. 
Similarly $\U_{\I,v}=\spec(\CC[z'_{-\alpha}]_{\alpha\in\Delta_{\I}^{+}})$ where again each $z'_{-\alpha}$
is an independent variable on which $\T$ acts via the weight $-\alpha$.  The $\T$-equivariant closed embedding
$\U_{\I,v}\hookrightarrow\U_{v}$ corresponds to a $\T$-equivariant surjective map of rings
$\CC[z_{-\alpha}]_{\alpha\in\invset_{v^{-1}}}\longrightarrow \CC[z'_{-\alpha}]_{\alpha\in\Delta_{\I}^{+}}$.  
If $\gamma$ is a weight in $\Span_{\ZZ_{\leq 0}}\I=\Span_{\ZZ_{\leq 0}}\Delta_{\I}^{+}$ 
then the dimension of the $\T$-eigenspace of weight $\gamma$ in both rings is the same.  
In particular, no monomial in the variables $\{z_{-\alpha}\}_{\alpha\in\Delta_{\I}^{+}}$ is in the kernel of the map,
while all monomials involving the variables $\{z_{-\alpha}\}_{\alpha\in\invset_{v^{-1}}\setminus\Delta_{\I}^{+}}$ are.
Therefore
the kernel of the surjection is the direct sum of the 
$\T$-eigenspaces consisting of the functions whose weight lies in 
$\SS$.
This proves ({\em b}).

\np
If $\lambda$ is dominant then $\L_{\lambda}$ is basepoint free on $\X$, and so the pullback map  $\psi_{v}^{*}$
from $\H^{0}(\X,\L_{\lambda})$ to $\H^{0}(\XI,\psi_{v}^{*}\L_{\lambda})=\H^{0}(\XI,\varphi_{v}^{*}(\L_{\lambda}|_{\X_v}))$
is nonzero.  On the other hand, by part ({\em a}) the pullback map $\psi_{v}^{*}$ is $\GI$-equivariant, and since
$\H^{0}(\XI,\psi_{v}^{*}\L_{\lambda})$ 
is an irreducible representation of $\GI$, 
$\psi_{v}^{*}$ must be surjective.  
The map $\varphi_{v}^{*}$ is therefore also surjective since $\psi_{v}^{*}$ factors through $\varphi_{v}^{*}$.
Under the composite map 
$\GI\longrightarrow\XI \stackrel{\varphi_{v}}{\longrightarrow}\X$ the point $1_{\GI}\in\GI$ gets sent to $v\in\X_{v}$, 
and hence the torus weight of $\varphi_{v}^{*}\L_{\lambda}$ at the image of $1_{\GI}$ in $\XI$ is $-v\lambda$, and 
therefore $\H^{0}(\XI,\varphi_{v}^{*}(\L_{\lambda}|_{\X_{v}})) \cong \V_{\mut}\otimes \LC_{\mup}$ as representations of $\ggI$ 
by Theorem \ref{thm:Borel-Weil-XI}, proving ({\em c}). \epf

\np
We will also need a variant of Proposition \ref{prop:XI-embedding}({\em a,c}) under the ``opposite'' hypothesis
that $\Delta_{\I}^{+}\cap\invset_{v^{-1}}=\emptyset$.  We omit the demonstration since it only involves minor 
modifications of the proof of Proposition \ref{prop:XI-embedding}.

\tpoint{Proposition} \label{prop:XI-embedding-opp}
Let $v$ be an element of the Weyl group such that 
that $\Delta_{\I}^{+}\cap\invset_{v^{-1}}=\emptyset$.  Then the map $\GI\longrightarrow\G$ defined by $g\mapsto gv$
induces a $\GI$-equivariant embedding $\psi'_{v}\colon\GI/\BI\longrightarrow \X$ sending $1_{\GI}$ to $v\in\X$.
For any dominant weight $\lambda$, 
$\H^{0}(\GI/\BI,\psi^{'*}_{v}\L_{\lambda}) \cong \V_{\mut}^{*}\otimes\LC_{\mup}$ as representations of $\ggI$, where
where $\mut$ and $\mup$ are the restrictions of $-v\lambda$ to $\gttil$ and $\ga$ respectively under the decomposition
$\gt=\gttil\oplus\ga$ from \S\ref{sec:GtilandGI}.

\np
The action of $\GI$ on $\XI$ factors through the center and so $\Gtil$ acts naturally on $\XI$.
As a $\Gtil$-variety $\XI$ (and $\GI/\BI$) are isomorphic in a unique way to $\Xtil:=\Gtil/\Btil$, and 
in the statement of the main theorem 
we will also use $\psi_{v}$ and $\psi'_{v}$ for the maps from $\Xtil$ 
into $\X$ given by the constructions in Propositions \ref{prop:XI-embedding} and \ref{prop:XI-embedding-opp}.

\np
{\bf Construction of $\G\times^{\B}\X_{\vv}$ and maps.}
For any sequence $\vv=(v_1,\ldots, v_{k+1})$ of Weyl group elements we set 
$\X_{\vv}:=\X_{v_1}\times\cdots\times\X_{v_{k+1}}$ and consider it as a $\B$-variety where $\B$ acts diagonally.
We define $\G\times^{\B}\X_{\vv}$ to be the quotient 
$\G\times^{\B}\X_{\vv}:=(\G\times\X_{\vv})/\B$ where the $\B$-action is given by 
$$b\cdot(g,x_1,\ldots,x_{k+1}) = (gb^{-1},b\cdot x_1,\ldots, b\cdot x_{k+1})$$
for a point $(g,x_1,\ldots, x_{k+1})$ of $\G\times\X_{v_1}\times\cdots\times\X_{v_{k+1}}$.

\np
The group $\G$ acts on $\G\times\X_{\vv}$ by left multiplication on the first factor.  Since this action commutes
with the action of $\B$ above it descends to an action of $\G$ on $\G\times^{\B}\X_{\vv}$.
The map from $\G\times\X_{\vv}$ to $\X^{k+1}$ given by 

\begin{equation}\label{eqn:Fv}
(g,x_1,\ldots, x_{k+1})\mapsto (g\cdot x_1,\ldots, g\cdot x_{k+1})
\end{equation}

\noindent
is invariant under the $\B$-action.  If we let $\G$ act on $\X^{k+1}$ diagonally then \eqref{eqn:Fv} is also 
$\G$-equivariant and hence descends to a $\G$-equivariant morphism 
$f_{\vv}\colon(\G\times^{\B}\X_{\vv})\longrightarrow\X^{k+1}$.

\np
Similarly, the map $\G\times\X_{\vv}\longrightarrow\G$ given by projection onto the first factor 
descends to a $\G$-equivariant map $\fo\colon(\G\times^{\B}\X_{\vv})\longrightarrow \X$ expressing 
$\G\times^{\B}\X_{\vv}$ as an $\X_{\vv}$-bundle over $\X$.  In particular, setting $\N=\dim(\X)=|\Delta^{+}|$, 
we obtain that $\dim(\X_{\vv})=\N+\sum_{i=1}^{k+1}\ell(v_i)$, 
and hence $\dim(\G\times^{\B}\X_{\vv})=\dim(\X^{k+1})$ if and only if $\sum_{i=1}^{k+1}\ell(v_i)=k\N$.

\tpoint{Proposition} \label{prop:degree} 
If $\vv=(v_1,\ldots, v_{k+1})$ and $\sum_{i=1}^{k+1} \ell(v_i) = k\N$ then the degree of 
$f_{\vv}\colon (\G\times^{\B}\X_{\vv}) \longrightarrow\X^{k+1}$ is given by the intersection number 
$\cap_{i=1}^{k+1}[\Omega_{\wo v^{-1}_i}] = \cap_{i=1}^{k+1}[\X_{v_i^{-1}}]$.

\bpf
After re-indexing $k$ as $k+1$, this is
\cite[Corollary (3.7.5)]{dr1}, along with the observation that the variety $\Q_{\vseq}$ used in the corollary
is our variety $\G\times^{\B}\X_{\vv}$, and the map $h\colon\Q_{\vseq}\longrightarrow\X^{k+1}$ considered there
is our map $f_{\vv}$.
\epf

\bpoint{The Littlewood-Richardson Cone}
\label{sec:LR-cone}
For any $k\geq 1$, let $\LRC(k)$ be the {\em Littlewood-Richardson cone}, i.e., 
the rational cone generated by the tuples $(\mu_1,\ldots, \mu_k,\mu)$ of dominant weights
such that $\V_{\mu}$ is a component of $\V_{\mu_1}\otimes\cdots\otimes\V_{\mu_k}$.  
It is known that $\LRC(k)$ is polyhedral.  
A face of $\LRC(k)$ is called {\em regular} if it intersects the locus of strictly dominant weights.

\np
{\bf Description of regular faces.}
For any set $\I$ of simple roots, we define $\P_{\I}$ to be the parabolic subgroup associated to $\I$. 
For any parabolic $\P\supseteq\B$ we denote the Weyl group of $\P$ by $\W_{\P}$.

\np
For a set $\I$ of simple roots we wish to consider elements $w_1$,\ldots, $w_k$, and $w$ 
of $\W$ satisfying the following conditions with respect to $\I$:

\begin{equation}\label{eqn:Iconditions}
\left\{
\mbox{
\begin{minipage}{0.8\textwidth}
\begin{itemize}
\item[({\em i})] Each $w_i$ is of minimal length in the coset $w_i\W_{\P_{\I}}$, and $w$ is of minimal length
in the coset $w\W_{\P_{\I}}$.
\item[({\em ii})] $\ell(w)=\sum_{i=1}^{k} \ell(w_i)$ and $\cap_{i=1}^{k}[\Omega_{w_i}]\cdot[\X_{w}]=1$.
\item[({\em iii})] 
The weight $\sum_{i=1}^{k} w_i^{-1}\cdot 0 - w^{-1}\cdot 0$ belongs to $\Span_{\ZZ_{\geq 0}} \I$,
\end{itemize}
\end{minipage}
}
\right.
\end{equation}

\np
where $u\cdot 0$ denotes the affine action of an element $u\in\W$ on the zero weight.
To produce examples of such $w_1$,\ldots, $w_k$, and $w$ it is usually easier to use the following equivalent
formulation of conditions \eqref{eqn:Iconditions}:

\begin{equation}\label{eqn:altIconditions}
\left\{
\mbox{
\begin{minipage}{0.8\textwidth}
\begin{itemize}
\item[({\em i})] The classes $[\Omega_{w_i}]$, $i=1$,\ldots, $k$,  and $[\Omega_{w}]$ are pullbacks of 
Schubert classes $\sigma_i$, $i=1$,\ldots, $k$ and $\sigma$ respectively from $\G/\PI$.
\item[({\em ii})] The coefficient of $\sigma$ when writing the product $\cap_{i=1}^{k}\sigma_i$ as a sum
of basis elements is $1$.
\item[({\em iii})] 
The weight $\sum_{i=1}^{k} w_i^{-1}\cdot 0 - w^{-1}\cdot 0$ belongs to $\Span_{\ZZ_{\geq 0}} \I$.
\end{itemize}
\end{minipage}
}
\right.
\end{equation}

\np
The conditions above are directly equivalent, i.e., \eqref{eqn:altIconditions}({\em i}) is equivalent to 
\eqref{eqn:Iconditions}({\em i}) and \eqref{eqn:altIconditions}({\em ii})  is equivalent to \eqref{eqn:Iconditions}({\em ii}).

\np
The work of Ressayre gives an explicit description of the regular faces of $\LRC(k)$.
The following is a translation of \cite[Theorem D]{Re} into our notation:

\tpoint{Theorem} \label{thm:Ressayre}

\begin{enumerate}
\item
Let $\I$ be a set of simple roots and 
$w_1$,\ldots, $w_k$, and $w$ elements of $\W$ satisfying conditions \eqref{eqn:Iconditions} with respect to $\I$.
Then the set 

$$
\left\{(\mu_1,\ldots,\mu_k,\mu)\in\LRC(k) \, | \, \sum_{i=1}^{k} w_i^{-1}\mu_i-w^{-1}\mu \in \Span_{\QQ_{\geq 0}} \I \right\}
$$

\np
is a regular face of codimension $(n-|\I|)$ of $\LRC(k)$.
Here $|\I|$ denotes the cardinality of the set $\I$ and $n$ the rank of $\G$.

\medskip
\item Any regular face of $\LRC(k)$ is of the form given in part ({\em a}).
\end{enumerate}

\np
The theorem of Ressayre above is not necessary for the proof of the reduction theorem.  Its importance for this
paper is that it links the combinatorial conditions used in the proof with the geometry of the 
Littlewood-Richardson cone, and guarantees that there are examples to which the reduction rules apply.

\section{Statement and proof of the reduction theorem}
\label{sec:proofofThm}

\setcounter{subsection}{1}
\setcounter{subsubsection}{0}

\np
\tpoint{Reduction Theorem} Suppose that we are given a set $\I$ of simple roots
 and elements $w_1$,\ldots, $w_k$, $w\in\W$ satisfying conditions  \eqref{eqn:Iconditions}({\em i},{\em ii}) 
 with respect to $\I$.
 Let $\Gtil$ be the semisimple part of $\PI$, $\Xtil=\Gtil/\Btil$, $\X=\G/\B$, and 
 $\psi:=\psi_{w_1^{-1}\wo}\times\cdots\times\psi_{w_k^{-1}\wo}\times\psi'_{w^{-1}}$ the $\Gtil$-equivariant map 
 $\Xtil^{k+1}\longrightarrow \X^{k+1}$ given by the constructions in 
 Propositions \ref{prop:XI-embedding} and \ref{prop:XI-embedding-opp}.
 Suppose that dominant weights $\mu_1$,\ldots, $\mu_k$, and $\mu$ satisfy

\begin{equation}\label{eqn:onface}
\sum_{i=1}^{k} w_i^{-1}\mu_i-w^{-1}\mu \in \Span_{\QQ} \I,
\end{equation}

\np
and
let $\mut_1$,\ldots, $\mut_k$, and $\mut$ be the reductions (cf.\ \S\ref{sec:GtilandGI})
of $w_1^{-1}\mu_1$, \ldots, $w_k^{-1}\mu_k$, and $w^{-1}\mu$ respectively to $\Ttil$. 
Set $\L:=\L_{-\wo\mu_1}\bt\cdots\bt\L_{-\wo\mu_k}\bt\L_{\mu}$ on $\X^{k+1}$. Then the
pullback of global sections of $\L$ by $\psi$ induces an isomorphism of vector spaces

$$(\V_{\smu_1}\otimes\cdots\otimes\V_{\smu_k}\otimes\V_{\smu}^{*})^{\G} \stackrel{\sim}{\longrightarrow}
(\V_{\smut_1}\otimes\cdots\otimes\V_{\smut_k}\otimes\V_{\smut}^{*})^{\Gtil},$$

\np
and in particular 
$\mult_{\G}(\V_{\smu},\V_{\smu_1}\otimes\cdots\otimes\V_{\smu_k}) = 
\mult_{\Gtil}(\V_{\smut},\V_{\smut_1}\otimes\cdots\otimes\V_{\smut_k})$.

\np
\bpf
We will construct a sequence of isomorphisms of vector spaces starting with 
$(\V_{\smu_1}\otimes\cdots\otimes\V_{\smu_k}\otimes\V_{\smu}^{*})^{\G}$ and ending with 
$(\V_{\smut_1}\otimes\cdots\otimes\V_{\smut_k}\otimes\V_{\smut}^{*})^{\Gtil}.$
Afterwards we will check that the composite isomorphism is that induced by pullback of global sections via $\psi$.

\np
{\bf Step 1.} 
Set $\lambda_i=-\wo\mu_i$ for $i=1$,\ldots, $k$, and $\lambda_{k+1}=\mu$. 
Let $\L$ be the line bundle $\L_{\lambda_1}\bt\cdots\bt\L_{\lambda_{k+1}}$ on $\X^{k+1}$ as above
so that $\H^{0}(\X^{k+1},\L) = \V_{\smu_1}\otimes\cdots\otimes\V_{\smu_k}\otimes\V_{\smu}^{*}$.

\np
Set $v_i=w_i^{-1}\wo$ for $i=1$,\ldots, $k$, $v_{k+1}=w^{-1}$, and $\vv=(v_1,\ldots, v_{k+1})$ and
consider the map $f_{\vv}\colon(\G\times^{\B}\X_{\vv})\longrightarrow\X^{k+1}$ from \S\ref{sec:Schub}. 
By Proposition 
\ref{prop:degree} the degree of $f_{\vv}$ is given by 
$$\bigcap_{i=1}^{k+1}[\Omega_{\wo v_i^{-1}}] = \bigcap_{i=1}^{k} [\Omega_{w_i}] \cdot [\X_{w}] 
\stackrel{\mbox{\scriptsize\eqref{eqn:Iconditions}({\em ii})}}{=} 1,$$

\np
and therefore $f_{\vv}$ is a proper birational 
map.  Since $\X^{k+1}$ is smooth it follows that $f_{\vv*}(f_{\vv}^{*}\L)=\L$ 
and therefore pullback induces an isomorphism
$$\H^{0}(\G\times^{\B}\X_{\vv},f_{\vv}^{*}\L)\stackrel{f_{\vv}^{*}}{\longleftarrow}\H^{0}(\X^{k+1},\L).$$

\np
Because $f_{\vv}$ is $\G$-equivariant, $f_{\vv}^{*}$ induces an isomorphism of $\G$-invariant subspaces, 
and we may therefore focus our attention on $\H^{0}(\G\times^{\B}\X_{\vv},f_{\vv}^{*}\L)^{\G}$.

\np
{\bf Step 2.}  Let $\Esh_{2}=\fopush (f_{\vv}^{*}\L)$, where $\fo\colon(\G\times^{\B}\X_{\vv})\longrightarrow\X$
is the map from \S\ref{sec:Schub}, 
and let $\E_2$ be the fibre of $\Esh_2$ over $e\in\X$. 
Since $\fo$ is also $\G$-equivariant, pushforward induces an isomorphism 
$\H^0(\G\times^{\B}\X_{\vv},f_{\vv}^{*}\L)^{\G}\stackrel{\sim}{\longrightarrow}\H^0(\X,\Esh_2)^{\G}$.
By Principle \ref{princip} restriction to the fibre over $e$ induces an isomorphism 
$\H^0(\X,\Esh_2)^{\G}  \stackrel{\sim}{\longrightarrow} \E_2^{\B}$.

\np
{\bf Step 3.} 
The fibre of $\fo$ over $e\in\X$ is $\X_{\vv}$.
Let $i_{\vv}\colon\X_{\vv}\longrightarrow\X^{k+1}$ be the restriction of $f_{\vv}$ to this fibre.
From the construction in \S\ref{sec:Schub} it follows that $i_{\vv}$ is the product of the natural 
inclusion maps 
$\X_{v_j}\hookrightarrow\X$ for $j=1$,\ldots, $k+1$. 
Hence by the theorem on cohomology and base change 
$\E_2=\H^0(\X_{\vv},i_{\vv}^{*}\L)$.
Set $\gamma=\sum_{i=1}^{k} w_i^{-1}\mu_i-w^{-1}\mu$; then 
$\gamma$ is the weight of $\T$ acting on $i_{\vv}^{*}\L$
at the point $p:=(v_1,\ldots,v_{k+1})\in\X_{\vv}$, and 
$\gamma\in\Span_{\QQ}\I$ by condition \eqref{eqn:onface}.

\np
Let $\U=\U_{v_1}\times\cdots\times\U_{v_{k+1}}$ be the product of the $\B$-stable affine spaces $\U_{v_i}$ from
\S\ref{sec:Schub}; the origin of $\U$ is the point $p$.
Since $\U$ is open in the irreducible variety $\X_{\vv}$, restriction
gives an $\B$-equivariant inclusion 
$\E_2=\H^{0}(\X_{\vv},i_{\vv}^{*}{\L}|_{\X_{\vv}}) \hookrightarrow \H^{0}(\U,i_{\vv}^{*}\L|_{\U}).$

\np
Set $\L_{\U}=(i_{\vv}^{*}\L)|_{\U}$.
Since $\U$ is isomorphic to affine space, $\L_{\U}$ is (non-equivariantly) trivial on $\U$. Let 
$\so$ be a section of $\L_{\U}$ which is nowhere vanishing. The torus $\T$ 
takes $\so$ to another nowhere vanishing
section which must therefore be a multiple of $\so$, i.e., $\T$ acts on $\so$ via a weight.  This must be the
same as the weight of the action of  $\L_{\U}$ at $p$, and so $\T$ acts on $\so$ with weight $\gamma$.
Let $\B^{+}$ be the unipotent radical of $\B$.  By the same reasoning,
$\B^{+}$ must take $\so$ to a multiple of itself. Since $\B^{+}$ has only the trivial 
one-dimensional representation $\so$ must be fixed by $\B^{+}$.

\np
Every section $s\in\H^{0}(\U,\L_{\U})$ can be written as $s=\so h$ for some function $h\in\H^{0}(\U,\Osh_{\U})$.
The section $s$ is $\B$-invariant if and only if $h$ is $\B^{+}$-invariant  and $h$ is an eigenfunction of $\T$
on which $\T$ acts via $-\gamma$.
For any weight $\delta$, let $\H^{0}(\U,\Osh_{\U})_{\delta}$ denote the space of 
eigenfunctions of $\T$ on which $\T$ acts via $\delta$.  
Let $\gb^{+}$ be the Lie algebra of $\B^{+}$ (i.e, the nilpotent radical of $\gb$); $\gb^{+}$ acts on 
$\H^{0}(\U,\Osh_{\U})$ via derivations.
By the correspondence above between 
sections of $\L_{\U}$ and functions on $\U$ we have 
$\H^0(\U,\L_{\U})^{\B}=\H^{0}(\U,\Osh_{\U})^{\gb^{+}}_{-\gamma}$.

\np
For each $\beta\in\delpos$ let $\die_{\beta}$ be a vector field giving the action of a nonzero element of
$\ggg^{\beta}\subseteq\gb^{+}$ 
on $\U$.  
Each $\die_{\beta}$ is a graded first-order differential operator of degree $\beta$, i.e., 
$$\die_{\beta}\left({\H^{0}(\U,\Osh_{\U})_{\delta}}\right)\subseteq\H^{0}(\U,\Osh_{\U})_{\delta+\beta}$$
for each weight $\delta$. 
Thus, we obtain 

\begin{equation}\label{eqn:H0B}
\H^{0}(\U,\L_{\U})^{\B}=
\H^{0}(\U,\Osh_{\U})_{-\gamma}^{\gb^{+}}  =
\bigcap_{\beta\in\delpos} \ker\left({
\H^{0}(\U,\Osh_{\U})_{-\gamma} \stackrel{\die_{\beta}}{\longrightarrow} \H^{0}(\U,\Osh_{\U})_{-\gamma+\beta}}\right).
\end{equation}

\np
By repeating the same argument with the subgroup $\BI$ we obtain a similar identification 

\begin{equation}\label{eqn:H0BI}
\H^{0}(\U,\L_{\U})^{\BI} =
\H^{0}(\U,\Osh_{\U})_{-\gamma}^{\gbI^{+}}  =
\bigcap_{\beta\in\delpos_{\I}} \ker\left({
\H^{0}(\U,\Osh_{\U})_{-\gamma} \stackrel{\die_{\beta}}{\longrightarrow} \H^{0}(\U,\Osh_{\U})_{-\gamma+\beta}}\right).
\end{equation}

\np
Since $\gamma\in\Span_{\QQ}{\I}$, if $\beta\in\Delta^{+}\setminus\Delta_{\I}^{+}$ 
then $-\gamma+\beta\not\in\Span_{\ZZ_{\leq 0}}\Delta^{+}$ and so $\H^{0}(\U,\Osh_{\U})_{-\gamma+\beta}=0$ by 
Lemma \ref{lem:zero-unless-positive}.
Thus the right-hand sides of \eqref{eqn:H0B} and \eqref{eqn:H0BI} are equal, and hence 
$\H^{0}(\U,\L_{\U})^{\B}=\H^{0}(\U,\L_{\U})^{\BI}$.  Since the inclusion map $\E_2\hookrightarrow\H^{0}(\U,\L_{\U})$
is $\B$-equivariant we conclude that $\E_2^{\B}=\E_2^{\BI}$.  Passing to 
the Lie algebra of $\BI$ we are reduced to studying $\E_2^{\gbI}$.

\np
{\bf Step 4.} 
An element $u\in\W$ is of minimal length in the coset $u\W_{\PI}$ if and only if $\Delta_{\I}^{+}\cap \invset_{u}=\emptyset$.  Applying this observation to each $w_i$, 
we conclude that 
$\Delta_{\I}^{+}\subseteq \invset_{w_i}^{\c}
=
\invset_{v_i^{-1}}$, and hence
by Proposition \ref{prop:XI-embedding}({\em b})
we have $\BI$-equivariant embeddings 
$\varphi_{v_i}\colon \XI\longrightarrow\X_{v_i}$ for $i=1$,\ldots, $k$.

\np
The variety $\X_{v_{k+1}}$ is stable under $\B$ and hence under the subgroup $\BI\subseteq\B$.
The stabilizer subgroup of $v_{k+1}\in\X$ is $v_{k+1}\B v_{k+1}^{-1}$ with roots 
$v_{k+1}\Delta^{+}=w^{-1}\Delta^{+}=-\invset_{w}\sqcup\invset_{w}^{\c}.$  Applying the observation on minimality
of length to $w$ we conclude that $\Delta_{\I}^{+}\subseteq\invset_{w}^{\c}$, 
and hence that $\BI\subseteq v_{k+1}\B v_{k+1}^{-1}$, i.e., $\BI$ fixes the point $v_{k+1}\in\X_{v_{k+1}}$.
Let $j_{v_{k+1}}\colon\spec(\CC)\longrightarrow\X_{v_{k+1}}$ be the $\BI$-equivariant inclusion of the point $v_{k+1}$.

\np
Finally, 
let $\varphi_{\vv}\colon \XI^{k}\longrightarrow \X_{\vv}$ be the map including $\XI^{k}=\XI^{k}\times\spec(\CC)$
into $\X_{\vv}$ via the product inclusions $\varphi_{v_1}\times\cdots\times\varphi_{v_{k}}\times j_{k+1}$ and set
$\E_3=\H^{0}(\XI^{k},\varphi_{\vv}^{*}i_{\vv}^{*}\L)$.  By the Kunnuth theorem 

\begin{equation}\label{eqn:Kunnuth}
\E_3=\left(\ott_{i=1}^{k}\H^0(\XI,\varphi_{v_i}^{*}(\L_{\lambda_i}|_{\X_{v_i}}))\right)\otimes 
\left(j_{k+1}^{*}(\L_{\lambda_{k+1}}|_{\X_{v_{k+1}}})\right).
\end{equation}

\np
By
Proposition \ref{prop:XI-embedding}({\em c}) each of the pullback maps 
$$\varphi_{v_i}^{*}\colon\H^0(\X_{v_i},\L_{\lambda_i}|_{\X_{v_i}})\longrightarrow
\H^0(\XI,\varphi_{v_i}^{*}(\L_{\lambda_i}|_{\X_{v_i}}))$$
are surjective for $i=1$,\ldots, $k$, 
and certainly $j_{k+1}^{*}\colon\H^0(\X_{v_{k+1}},\L_{\lambda_{k+1}}|_{\X_{v_{k+1}}})\longrightarrow 
j_{k+1}^{*}(\L_{\lambda_{k+1}}|_{v_{k+1}})$ is 
surjective since $\L_{\lambda_{k+1}}$ is basepoint free on $\X$ and the pullback is to a point.
Thus the $\BI$-equivariant pullback map $\varphi_{\vv}^{*}\colon\E_2\longrightarrow\E_3$ is
surjective.   We want to see that this surjection induces an isomorphism of $\gbI$-invariants.

\np
Let $\E_1$ be the kernel of the surjection above.  If $\Ish$ is the ideal sheaf of $\varphi_{\vv}(\XI^{k})$ in $\X_{\vv}$
then $\E_1=\H^{0}(\X_{\vv},(i_{\vv}^{*}\L)\otimes_{\Osh_{\X_{\vv}}}\Ish)$. As in step 3 we will
analyze $\E_1$ via the inclusion $\E_1\hookrightarrow\H^{0}(\U,\L_{\U}\otimes_{\Osh_{\U}}\Ish|_{\U})$ obtained
by restriction to $\U$.
As in step 3 every section
$s\in\H^{0}(\U,\L_{\U}\otimes_{\Osh_{\U}}\Ish|_{\U})$ can be written as $\so h$ with $h\in\H^{0}(\U,\Ish|_{\U})$.
Since $\XI^{k}$ is a product subvariety in the product variety $\X_{\vv}$, and $\U$ is a product subset, the 
ideal $\H^{0}(\U,\Ish|_{\U})$ is the sum of the pullbacks to $\U$ of the 
ideals of $\XI|_{\U_{v_i}}$, $i=1$,\ldots, $k$ and the ideal 
of the point $v_{k+1}\in\U_{v_{k+1}}$.  By Proposition
\ref{prop:XI-embedding}({\em b}) for each $i=1$,\ldots, $k$, 
the ideal of $\XI|_{\U_{v_i}}$ consists of the direct sum of the $\T$-eigenspaces of functions on $\U_{v_i}$ with
torus weights contained in $\SS=\left(\Span_{\ZZ_{\leq0}}(\Delta^{+}\setminus\Delta_{\I}^{+})\right)\setminus\{0\}$.   
Now
$\U_{v_{k+1}}=\spec(\CC[z_{-\alpha}]_{\alpha\in\invset_{v_{k+1}^{-1}}})$ and the ideal of $v_{k+1}$ in $\U_{v_{k+1}}$ 
is generated by 
$\{z_{-\alpha}\}_{\alpha\in\invset_{v_{k+1}^{-1}}}$.  Since $\invset_{v_{k+1}^{-1}}=\invset_{w}$, and again using
the observation on the minimality of $w$, we conclude that the weights of 
all $\T$-eigenfunctions in the ideal of $v_{k+1}$ in 
$\U_{v_{k+1}}$ are also contained in $\SS$. 
Pulling back these ideals to $\U$, and using the fact that 
$\T$ acts on $\so$ with
weight $\gamma\in\Span_{\QQ}\I$, we conclude that all $\T$-eigensections 
$s=\so h\in\H^{0}(\U,\L_{\U}\otimes_{\Osh_{\U}}\Ish|_{\U})$ have weights outside $\Span_{\ZZ}\I$. 
Since 
$\E_1\hookrightarrow \H^{0}(\U,\L_{\U}\otimes_{\Osh_{\U}}\Ish|_{\U})$  is a $\BI$-equivariant inclusion, we
conclude that the same is true for the weights of $\E_1$. In particular, 
no weight of $\E_1$ is contained in $\{0\}\cup\Delta_{\I}^{+}$.  Thus by Lemma 
\ref{lem:bI-invts} the surjection $\E_2\longrightarrow\E_3$ 
induces an isomorphism $\E_2^{\gbI}\stackrel{\sim}{\longrightarrow}\E_3^{\gbI}$.

\np
{\bf Step 5.}
By Lemma \ref{prop:XI-embedding}({\em c}), for $i=1$,\ldots, $k$ we have 
$$\H^0(\XI,\varphi_{v_i}^{*}(\L_{\lambda_i}|_{\X_{v_i}}))\cong\V_{\mut_i}\otimes\LC_{\mup_i}$$
as $\ggI$-modules
and the $\gbI$-module structure on 
$\H^0(\XI,\varphi_{v_i}^{*}(\L_{\lambda_i}|_{\X_{v_i}}))$
is simply the restriction of the $\ggI$-module structure.
Here $\mut_i$ and $\mup_i$ are  restrictions to $\gttil$ and $\ga$ respectively of $w_i^{-1}\mu_i= -v_i\lambda_i$
using the decomposition $\gt=\gttil\oplus\ga$ from \S\ref{sec:GtilandGI}.

\np
Similarly, the one-dimensional $\gt$-representation $j_{k+1}^{*}(\L_{\lambda_{k+1}}|_{\X_{v_{k+1}}})$ decomposes
as $\LC_{-\mut}\otimes\LC_{-\mup}$ where $\mut$ and $\mup$ are the restrictions to $\gttil$ and $\ga$ respectively
of $v_{k+1}\lambda_{k+1}=w^{-1}\mu$.

\np
Thus, using \eqref{eqn:Kunnuth} and collecting the one-dimensional representations of $\ga$ we have 

$$\E_3 = 
\left(\ott_{i=1}^{k} \V_{\mut_i}\right)\otimes\LC_{-\mut}\otimes\LC_{\left(\sum_{i=1}^{k}\mup_i\right) - \mup}.
$$

\np
However, since restriction is a homomorphism, $\left(\sum_{i=1}^{k} \mup_i\right)-\mup$ is just the restriction to $\ga$
of the weight $\gamma$, and this is zero by Condition \eqref{eqn:onface} and Lemma \ref{lem:cancelling}.
Thus 

\begin{equation}\label{eqn:E3decomp}
\E_3 = 
\left(\ott_{i=1}^{k} \V_{\mut_i}\right)\otimes\LC_{-\mut}
\end{equation}

\np
and hence $\E_3$ is a $\gbI$-module with trivial $\ga$-action, i.e., 
$\E_3$
is really a $\gbI/\ga=\gbtil$-module 
and so $\E_3^{\gbI} = \E_3^{\gbtil}$.

\np
{\bf Step 6.} It is straightforward to see that 
$\E_3^{\gbtil} = (\V_{\mut_i}\otimes\cdots\otimes\V_{\mut_k}\otimes\V_{\mut}^{*})^{\Gtil}$ which will finish the
construction of the isomorphism.

\np
The most direct argument is to notice that the $\gbtilp$-invariants of $\E_3$ are, by 
\eqref{eqn:E3decomp}, the highest-weight subspaces of the irreducible components of $\otimes_{i=1}^{k}\V_{\mut_i}$
tensored with $\LC_{-\mut}$, and
hence the $\gbtil$-invariants of $\E_3$ are the subspace of highest-weight vectors of weight $\mut$ in
$\otimes_{i=1}^{k}\V_{\mut_i}$, which is precisely the subspace  
$(\V_{\mut_i}\otimes\cdots\otimes\V_{\mut_k}\otimes\V_{\mut}^{*})^{\Gtil}$.

\np
A more geometric approach, inducing the isomorphism of vector spaces directly, is to let $\Esh_3$ be the vector bundle
on $\Xtil=\Gtil/\Btil$ whose fibre over $e\in\Xtil$ is $\E_3$.  By Principle \ref{princip} 
$\E_3^{\gbtil}=\E_3^{\Btil}=\H^{0}(\Xtil,\Esh_3)^{\Gtil}$.  
Equation \eqref{eqn:E3decomp} shows that $\Esh_3=(\ott_{i=1}^{k}\V_{\mut_i})\otimes_{\Osh_{\Xtil}}\L_{\mut}$, and
hence 

$$\H^{0}(\Xtil,\Esh_3)= \left(\ott_{i=1}^{k}\V_{\mut_i}\right)\otimes\H^{0}(\Xtil,\L_{\mut})  =
\V_{\mut_i}\otimes\cdots\otimes\V_{\mut_k}\otimes\V_{\mut}^{*}
$$

\np
by the Borel-Weil theorem. Taking $\Gtil$-invariants finishes the alternate argument for Step 6 and the construction
of the isomorphism.

\np
{\bf Composition of steps 1--6.}
Finally, we want to check that the composite isomorphism is that induced by pullback via $\psi$.  
Recall that we are identifying $\Xtil$ and $\XI$ by the unique isomorphism respecting their structure 
as $\Gtil$-varieties.
Let $\Ltil=\psi^{*}\L$.  It is straightforward to check (c.f.\ Propositions 
\ref{prop:XI-embedding} and \ref{prop:XI-embedding-opp}) that 
$\H^{0}(\Xtil^{k+1},\Ltil)=\V_{\mut_1}\otimes\cdots\otimes\V_{\mut_k}\otimes\V_{\mut}^{*}$.
Let $\pitil\colon\Xtil^{k+1}\longrightarrow\Xtil$ be projection onto the final factor.  
Pushing
forward by $\pitil$ we obtain 
$\H^{0}(\Xtil^{k+1},\Ltil)=\H^{0}(\Xtil,\pitil_{*}\Ltil)$.  The main point is that
$\pitil_{*}\Ltil=\Esh_3$ and that the pullback map $\psi^{*}$ on global sections induces
the isomorphism $\H^{0}(\X^{k+1},\L)^{\G}\stackrel{\sim}{\longrightarrow}\H^{0}(\X,\Esh_3)^{\Gtil}$ obtained by
combining steps 1 through 6.

\np
To see this, let $\Xtil^{k}$ be the fibre of $\pitil$ over $\etil\in\Xtil$, and 
let $j\colon \Xtil^{k}\longrightarrow\Xtil^{k+1}$ be the inclusion of this fibre in $\Xtil^{k+1}$.
By the theorem on cohomology and base change, the fibre of $\pitil_{*}\Ltil$ over $\etil$ is equal to 
$\H^{0}(\Xtil^{k},\Ltil|_{\Xtil^{k}})$.  By a straightforward check the composite map $\psi\circ j$ is equal to 
$i_{\vv}\circ\varphi_{\vv}$ and hence 
$\H^{0}(\Xtil^{k},\Ltil|_{\Xtil^{k}})=
\H^{0}(\Xtil^{k},\varphi_{\vv}^{*}i_{\vv}^{*}\L)=\E_3$ by the definition in step 4. Thus $\pi_{*}\Ltil=\Esh_{3}$.
The content of steps 1--5 is that restriction to $\Xtil^{k}$ (i.e., pullback by $\psi\circ j$)
induces an isomorphism $\H^{0}(\X^{k+1},\L)^{\G}\cong\E_3^{\Btil}$.
Since $\psi$ is $\Gtil$-equivariant, $\G$-invariant sections pull back to $\Gtil$-invariant sections, and so 
the composite isomorphism from steps 1--5 factors as

$$
\H^{0}(\X^{k+1},\L)^{\G} \stackrel{\psi^{*}}{\longrightarrow}\H^0(\Xtil^{k+1},\Ltil)^{\Gtil}
\stackrel{j^{*}}{\longrightarrow} \H^{0}(\Xtil^{k},\Ltil|_{\Xtil^{k}})^{\Btil}
=\E_3^{\Btil}.
$$

\np
Via the identification $\H^{0}(\Xtil^{k+1},\Ltil)^{\Gtil}=\H^{0}(\Xtil,\Esh_3)^{\Gtil}$ the map induced by 
$j^{*}$ is simply the natural restriction $\H^{0}(\Xtil,\Esh_3)^{\Gtil}\longrightarrow\E_3^{\Btil}$, which is
an isomorphism by Principle \ref{princip}.
The isomorphism $\E_3^{\Btil}\cong \H^{0}(\Xtil,\Esh_3)^{\Gtil}$ in step 6 is simply its
inverse.  Thus the map $\H^{0}(\X^{k+1},\L)^{\G}\longrightarrow\H^{0}(\Xtil^{k+1},\Ltil)^{\Gtil}$ induced by
pullback by $\psi$ is the composition of the maps from steps 1--6, and in particular is an isomorphism.
This finishes the proof of the reduction theorem. 
\epf

\np
{\bf Remarks.} Note that $w_1$,\ldots, $w_k$, and $w$ do not have to satisfy \eqref{eqn:Iconditions}({\em iii}) 
in order to apply the reduction theorem. 
Without \eqref{eqn:Iconditions}({\em iii}) 
however it is not clear that there are examples where the reduction rule applies, whereas such examples
are guaranteed by Theorem \ref{thm:Ressayre} if all the conditions do hold.  In applications of the reduction theorem,
it is convenient that one only has to verify the condition 
$\sum_{i=1}^{k} w_i^{-1}\mu_i-w^{-1}\mu \in \Span_{\QQ} \I$ and not that the sum is in $\Span_{\QQ_{\geq 0}}\I$.

\tpoint{Corollary} Suppose that $w_1$,\ldots, $w_k$, and $w$ satisfy \eqref{eqn:Iconditions}({\em i}) with respect
to $\I$, and that 
$\cap_{i=1}^{k}[\Omega_{w_i}]\cdot[\X_{w}]\neq 0$ (i.e., instead of $=1$).  
Then for any dominant weights $\mu_1$,\ldots, $\mu_k$, and $\mu$ such that 
$$
\sum_{i=1}^{k} w_i^{-1}\mu_i-w^{-1}\mu \in \Span_{\QQ} \I,$$

\np
we have 
$\mult_{\G}(\V_{\smu},\V_{\smu_1}\otimes\cdots\otimes\V_{\smu_k})  \leq
\mult_{\Gtil}(\V_{\smut},\V_{\smut_1}\otimes\cdots\otimes\V_{\smut_k})$, where $\mut_1$,\ldots, $\mut_k$, and
$\mut$ are the restrictions to $\Ttil$ of $w_1^{-1}\mu_1$,\ldots, $w_k^{-1}\mu_k$ and $w^{-1}\mu$ respectively.

\bpf
We repeat the proof of the reduction theorem.  The only difference occurs in Step 1, since the map $f_{\vv}$ now
may  have degree greater than one, and so we can only conclude that $f_{\vv}^{*}$ induces an 
inclusion 
$$\H^{0}(\G\times^{\B}\X_{\vv},f_{\vv}^{*}\L)^{\G}\stackrel{f_{\vv}^{*}}{\hookleftarrow}\H^{0}(\X^{k+1},\L)^{\G}
\stackrel{\sim}{\longleftarrow}(\V_{\smu_1}\otimes\cdots\otimes\V_{\smu_k}\otimes\V_{\smu}^{*})^{\G}.$$

\np
Following through the rest of the steps, we obtain an isomorphism 
$$\H^{0}(\G\times^{\B}\X_{\vv},f_{\vv}^{*}\L)^{\G}\stackrel{\sim}{\longrightarrow}
(\V_{\smut_1}\otimes\cdots\otimes\V_{\smut_k}\otimes\V_{\smut}^{*})^{\Gtil}
$$
and taking dimensions gives the inequality. \epf

\section{Examples}
\label{sec:examples}

\point
In this section we work out a number of explicit examples of reduction rules.
The rules in \S\ref{sec:Sl6}--\ref{sec:GHreduction} are of type $\A$, and so already covered by the results in
\cite{dw} and \cite{ktt} (the rule in \S\ref{sec:GHreduction} actually predates those papers -- it is
due to Griffiths and Harris).  
However the notation used in those papers is different from ours (the rules are expressed in $\Gl_{n}$ weights, 
and the combinatorial data describing the regular faces is presented in a different form) and the examples 
are included partly to compare the two approaches.  

\np
To check if a weight is in $\Span_{\QQ}\I$ one simply converts from the basis of fundamental weights to the root
basis by multiplying by the inverse transpose of the Cartan matrix, and then checks that the coordinates of all
simple roots outside of $\I$ are zero.  This is mentioned again in the first example, but afterwards we just  
write out the corresponding condition.

\np
In order to check that \eqref{eqn:Iconditions}({\em iii}) holds, the formula

\begin{equation}\label{eqn:rootsum}
w^{-1}\cdot 0 
=
w^{-1}\rho-\rho
= -\sum_{\alpha\in\invset_{w}} \alpha
\end{equation}

\np
is useful. Mostly, however we will also omit the explicit calculation checking this condition.  In particular, in 
type $\A_n$ when $|\I|=n-1$, condition \eqref{eqn:Iconditions}({\em iii}) follows from the condition 
$\sum_i \ell(w_i) = \ell(w)$ in \eqref{eqn:Iconditions}({\em ii}), and so does not need to be checked again.

\np
Because the reduction rules (and the multiplicities) depend only on the type of the group, we will label 
the examples and $\mult$ by the type, 
the only exception being for examples involving $\Gl_{n+1}$.
The labelling of the roots 
follows the usual convention in \cite[Chapter VI]{Bo}.
We will use $\alpha_1$,\ldots, $\alpha_n$ for the simple roots, and $s_1$,\ldots, $s_n$ for the corresponding
simple reflections.
After each of the examples we give an explicit instance with strictly dominant weights where the rule applies. 
By Theorem \ref{thm:Ressayre} such instances always exist.

\bpoint{An $\A_5$ to $\A_2\times\A_2$ reduction rule}\label{sec:Sl6}
Let $\G$ be of type $\A_5$ and $\I=\{\alpha_1,\alpha_2,\alpha_4,\alpha_5\}$ so 
that $\G/\PI=\Gr(2,5)$, the Grassmanian of two-planes in $\PP^5$. 
The Schubert basis for $\H^{*}(\Gr(2,5),\ZZ)$ consists of the classes $\sigma_{a_1,a_2,a_3}$ with 
$2\geq a_1\geq a_2\geq a_3\geq 0$.
In $\H^{*}(\Gr(2,5),\ZZ)$ we have the well-known cohomology multiplication $\sigma_{1,0,0}\cdot\sigma_{1,0,0}=
\sigma_{2,0,0}+\sigma_{1,1,0}$.
The pullback of $\sigma_{1,0,0}$ to $\X=\G/\B$ is $[\Omega_{s_3}]$ and the pullback of $\sigma_{2,0,0}$ to
$\X$ is $[\Omega_{s_4s_3}]$, so that if we pick $w_1=w_2=s_3$ and $w=s_4s_3$ then 
$w_1$, $w_2$, and $w$ satisfy \eqref{eqn:Iconditions} with respect to $\I$.  The group $\Gtil$ we are reducing to is
of type $\A_2\times \A_2$, obtained by deleting the middle node of the Dynkin diagram for $\G$.

\np
If $\mu_1=(a_1,a_2,a_3,a_4,a_5)$, $\mu_2=(b_1,b_2,b_3,b_4,b_5)$, and 
$\mu=(c_1,c_2,c_3,c_4,c_5)$ then 

\begin{eqnarray*}
w_1^{-1}\mu_1& = &  (a_1, a_2+a_3, -a_3, a_3+a_4, a_5) \\
w_2^{-1}\mu_2 & = &  (b_1, b_2+b_3, -b_3, b_3+b_4, b_5) \\
w^{-1}\mu  & = &   (c_1,c_2+c_3+c_4,-c_3-c_4,c_3,c_4+c_5)  \\
\end{eqnarray*}

\np
The group $\Gtil$ is a product group and we will use ``$|$'' to indicate the division of the restricted weight
among the two factors.  Since we are deleting the middle node of the Dynkin diagram, the restriction is
obtained by ignoring the middle coefficients in the formulas above, so that
$\mut_1=(a_1, a_2+a_3 \st a_3+a_4, a_5)$, $\mut_2=(b_1,b_2+b_3 \st b_3+b_4,b_5)$, and
$\mut=(c_1,c_2+c_3+c_4 \st c_3,c_4+c_5)$.

\np
The condition that the point $(\mu_1,\mu_2,\mu)$ lie on the face of $\LRC(2)$ determined by $\I$ and $w_1$, $w_2$, 
and $w_3$ is that the coefficient of $\alpha_3$ is zero when writing $w_1^{-1}\mu_1+w_2^{-1}\mu_2-w^{-1}\mu$ as
a sum of simple roots (with $\QQ$-coefficients).   This is easily computed by multiplying the sum, 
in the coordinates of the fundamental weights as above, by the inverse transpose of the Cartan matrix for $\A_5$ and 
looking at the middle coefficient.  This coefficient is 

$$ \frac{1}{2}\left({(a_1 + 2 a_2  + a_3 + 2a_4 + a_5) + (b_1+2b_2+b_3+2b_4+b_4) - (c_1+2c_2+c_3+c_5)\rule{0cm}{0.4cm}}\right),$$

\np
and thus we arrive at our first example of a reduction rule. 

\np
{\bf Reduction rule:} If

$$ c_1+2c_2+c_3+c_5 = (a_1 + 2 a_2  + a_3 + 2a_4 + a_5) + (b_1+2b_2+b_3+2b_4+b_4)$$

\np
then
$\displaystyle \mult_{\A_5}(\V_{\mu},\V_{\mu_1}\otimes\V_{\mu_2}) 
= \mult_{\A_2\times\A_2}(\V_{\mut},\V_{\mut_1}\otimes\V_{\mut_2})$, where $\mut_1$, $\mut_2$, and $\mut$ are 
given by the formulas above.

\np
{\bf Example:} $\mu_1=(4,2,10,6,10)$, $\mu_2=(10,4,12,4,2)$, $\mu=(10,22,1,1,25)$, $\mut_1=(4,12\st 16,10)$, 
$\mut_2=(10,16\st 16,2)$, $\mut=(10,24\st 1,26)$;  the multiplicity is $10$.

\np
In $\Gl_6$ weights, the rule has the following form.

\np
{\bf Reduction rule:} 
If dominant $\Gl_6$ weights 
$\mu_1=(a_0,a_1,a_2,a_3,a_4,a_5)$,$\mu_2=(b_0,b_1,b_2,b_3,b_4,b_5)$ and $\mu=(c_0,c_1,c_2,c_3,c_4,c_5)$ (which
we assume satisfy $\sum_{i} c_i = \sum_{i} a_i + \sum_{i} b_i$) also satisfy

$$c_0+c_1+c_4 = (a_0+a_1+a_3) + (b_0+b_1+b_3)$$

\np
then
$\displaystyle \mult_{\Gl_6}(\V_{\mu},\V_{\mu_1}\otimes\V_{\mu_2}) 
= \mult_{\Gl_3\times\Gl_3}(\V_{\mut},\V_{\mut_1}\otimes\V_{\mut_2})$
where 
$\mut_1=(a_0,a_1,a_3\st a_2,a_4,a_5)$, $\mut_2=(b_0,b_1,b_3 \st b_2,b_4,b_5)$, and 
$\mut=(c_0,c_1,c_4 \st c_2,c_3,c_5)$.

\np
{\bf Example:} $\mu_1=(32,28,26,16,10,0)$, $\mu_2=(32,22,18,6,2,0)$, $\mu=(60,51,28,26,25,2)$, 
$\mut_1=(32,28,16\st 26,10,0)$, $\mut_2=(32,22,6\st 18,2,0)$, $\mut=(60,51,25\st 28,26,2)$; multiplicity is 12.

\bpoint{An $\A_{n}$ to $\A_{n-1}$ reduction rule}  \label{sec:GHreduction}
Let $\G$ be of type $\A_n$ and $\I=\{\alpha_2,\alpha_3,\ldots, \alpha_n\}$ so that $\G/\PI=\PP^{n}$. 
We have $\H^{*}(\PP^{n},\ZZ)=\ZZ[h]/(h^{n+1})$, where $h\in\H^{2}(\PP^{n},\ZZ)$ is the hyperplane class.
Each $h^{i}$ ($1\leq i \leq n$) pulls back to the class $[\Omega_{s_is_{i-1}\cdots s_{1}}]$ in the cohomology ring of $\X$.  For any $0\leq i,j,k\leq n$ with $i+j=k$ we have the obvious cohomology multiplication $h^{i}\cdot h^{j}=h^{k}$.
Setting $w_1=s_is_{i-1}\cdots s_1$, $w_2=s_{j}s_{j-1}\cdots s_1$, and $w=s_ks_{k-1}\cdots s_1$, then $w_1$, $w_2$, and
$w$ satisfy \eqref{eqn:Iconditions} with respect to $\I$.  The group $\Gtil$ we are reducing to is of type $\A_{n-1}$
obtained by deleting the first node in the Dynkin diagram for $\G$.

\np
If $\mu_1=(a_1,\cdots, a_n)$, $\mu_2=(b_1,\ldots, b_n)$ and $\mu=(c_1,\cdots, c_n)$ are dominant weights then 

\begin{eqnarray*}
w_1^{-1}\mu_1& = &  (-a_1-a_2-\cdots-a_i,a_1,a_2,\cdots, a_{i-1},a_i+a_{i+1},a_{i+2},\cdots, a_n), \\
w_2^{-1}\mu_2& = &  (-b_1-b_2-\cdots-b_j,b_1,b_2,\cdots, b_{j-1},b_j+b_{j+1},b_{j+2},\cdots, b_n), \\
w^{-1}\mu& = &  (-c_1-c_2-\cdots-c_k,c_1,c_2,\cdots, c_{k-1},c_k+c_{k+1},c_{k+2},\cdots, c_n). \\
\end{eqnarray*}

\np
Restriction to $\Gtil$ simply ignores the first entries, so 

\begin{equation}\label{eqn:GHreductions}
\left\{
\mbox{
\begin{minipage}{0.5\textwidth}
\begin{eqnarray*}
\mut_1 & = & (a_1,\cdots, a_{i-1},a_{i}+a_{i+1},a_{i+2},\cdots, a_n), \\
\mut_2 & = & (b_1,\cdots, b_{j-1},b_{j}+b_{j+1},b_{j+2},\cdots, b_n), \,\,\mbox{and}\\
\mut & = & (c_1,\cdots, c_{k-1},c_{k}+c_{k+1},c_{k+2},\cdots, c_n).
\end{eqnarray*}
\end{minipage}
}
\right.
\end{equation}

\np
Here (and above) coefficients with indices greater than $n$ are assumed to be zero.

\np
Writing 
$w_1^{-1}\mu_1+w_2^{-1}\mu_2-w^{-1}\mu$ as a sum of simple roots and multiplying by $n+1$ to clear denominators, 
the coefficient of $\alpha_1$ is 

$$
(n+1)\sum_{r=i+1}^{n} a_r - \sum_{r=1}^{n} r a_r +
(n+1)\sum_{r=j+1}^{n} b_r - \sum_{r=1}^{n} r b_r -
(n+1)\sum_{r=k+1}^{n} c_r + \sum_{r=1}^{n} r c_r.
$$

\np
Thus we obtain the following family of reduction rules.

\np
{\bf Reduction rule:} For any integers $0\leq i,j,k\leq n$ with $i+j=k$, if 
dominant weights $\mu_1=(a_1,\cdots, a_n)$, $\mu_2=(b_1,\ldots, b_n)$ and $\mu=(c_1,\cdots, c_n)$ satisfy

\begin{equation}\label{eqn:GHconditions}
\rule{0.5cm}{0cm}
(n+1)\sum_{r=k+1}^{n} c_r - \sum_{r=1}^{n} r c_r =
(n+1)\sum_{r=i+1}^{n} a_r - \sum_{r=1}^{n} r a_r +
(n+1)\sum_{r=j+1}^{n} b_r - \sum_{r=1}^{n} r b_r 
\end{equation}

\np
then 
$\mult_{\A_{n}}(\V_{\smu},\V_{\smu_1}\otimes\V_{\smu_2})=\mult_{\A_{n-1}}(\V_{\smut},\V_{\smut_1}\otimes
\V_{\smut_2})$, where $\mut_1$, $\mut_2$, and $\mut$ are given by \eqref{eqn:GHreductions}.

\np
{\bf Example:} $n=5$, $i=j=1$, $k=2$, $\mu_1=(3,1,3,2,1)$, $\mu_2=(4,1,2,3,4)$, $\mu=(1,1,8,3,4)$,
$\mut_1=(4,3,2,1)$, $\mut_2=(5,2,3,4)$, $\mut=(1,9,3,4)$; the multiplicity is 24.

\np
This rule is much cleaner in $\Gl_{n+1}$ coordinates.  

\np
{\bf Reduction rule:} If 
$\mu_1=(a_0,\ldots, a_n)$, $\mu_2=(b_0,\ldots, b_n)$, and $\mu=(c_0,\cdots, c_n)$ are dominant $\Gl_{n+1}$ weights 
(again with $\sum c_i = \sum a_i + \sum b_i$), and $0\leq i,j,k\leq n$ such that $i+j=k$, then if $c_k=a_i+b_j$ we
have $\mult_{\Gl_{n+1}}(\V_{\smu},\V_{\smu_1}\otimes\V_{\smu_2}) = 
\mult_{\Gl_{n}}(\V_{\smut},\V_{\smut_1}\otimes\V_{\smut_2})$, where
$\mut_1$, $\mut_2$, and $\mut$ are obtained by deleting the entries $a_i$, $b_j$, and $c_k$ from $\mu_1$, $\mu_2$,
and $\mu$ respectively. 

\np
{\bf Example:} $n=6$, $i=1$, $j=2$, $k=3$, 
$\mu_1=(16,13,12,9,7,3,0)$, $\mu_2=(21,16,13,12,9,5,0)$, $\mu=(29,28,27,26,13,9,4)$,
$\mut_1=(16,12,9,7,3,0)$, $\mut_2=(21,16,12,9,5,0)$, $\mut=(29,28,27,13,9,$ $4)$; the multiplicity is 108.

\np
This $\Gl_{n+1}$ rule appears as Reduction Formula I for Schubert calculus in \cite[p.\ 202]{GH}.
(The rule given there
does not appear exactly as stated above, but is equivalent to it after making the translation from 
intersecting three Schubert cycles to computing the multiplicity of a representation in a tensor product, and after
using the indexing for the fundamental weights starting with zero.)

\bpoint{A three-factor reduction rule} The most important case for Littlewood-Richardson problems (i.e., the problem
of computing $\mult_{\G}(\V_{\smu},\V_{\smu_1}\otimes\cdots\otimes\V_{\smu_k})$) is the case with two factors, as
in the examples above.  The main theorem, however, gives the construction of reduction rules for an arbitrary number
of factors, and we give a three-factor example here.  For simplicity, we just repeat the $\Gl_{n+1}$ to $\Gl_{n}$
reduction in \S\ref{sec:GHreduction}, but now using the multiplication 
$h^{i}\cdot h^{j}\cdot h^{k}=h^{m}$ in $\H^{*}(\PP^{n},\ZZ)$ whenever $0\leq i,j,k,m\leq n$ and $m=i+j+k$.  
This gives:

\np
{\bf Reduction rule:}
For any $0\leq i,j,k,m\leq n$ with $i+j+k=m$, then for any dominant $\GL_{n+1}$ weights
$\mu_1=(a_0,\ldots, a_n)$, $\mu_2=(b_0,\ldots, b_n)$, $\mu_3=(c_0,\ldots, c_n)$, and $\mu=(d_0,\ldots, d_n)$, 
if $d_m=a_i+b_j+c_k$ then 
$\mult_{\Gl_{n+1}}(\V_{\smu},\V_{\smu_1}\otimes\V_{\smu_2}\otimes\V_{\smu_3}) = 
\mult_{\Gl_{n}}(\V_{\smut},\V_{\smut_1}\otimes\V_{\smut_2}\otimes\V_{\smut_3})$, where 
$\mut_1$, $\mut_2$, $\mut_3$, and $\mut$ are obtained by deleting the entries $a_i$, $b_j$, $c_k$ and $d_m$ from 
$\mu_1$, $\mu_2$, $\mu_3$, and $\mu$ respectively.  This rule generalizes to a larger number of factors in the
obvious way. 

\np
{\bf Example:} $n=4$, $i=j=k=1$, $m=3$, 
$\mu_1=(36,28,24,16,0)$, 
$\mu_2=(40,24,20,8,0)$,
$\mu_3=(94,14,11,9,0)$,
$\mu=(118,68,67,66,5)$,
$\mut_1=(36,24,16,0)$,
$\mut_2=(40,20,8,0)$,
$\mut_3=(94,11,9,0)$,
$\mut=(118,68,67,5)$;
the multiplicity is 196.














\np
Even though the Littlewood-Richardson coefficients for the decomposition of the tensor product of two
irreducible representations determine the coefficients for the decomposition of the tensor product of 
$k$ irreducible representation, there does not seem to be an obvious argument for 
deducing the $k$-factor reduction rules from the two-factor reduction rules.

\bpoint{A codimension-two reduction}  \label{sec:codim-two}
The previous examples have all been codimension-one reductions, i.e.,
starting with a codimension-one regular face of $\LRC(k)$ we obtain a rule with a single condition to 
check which reduces the rank of the group by one.  In this section we give a codimension-two example.  
By Corollary \ref{cor:codimr-restriction} below, 
any codimension-$r$ rule can be obtained as a succession of $r$ codimension-one rules, 
but it is sometimes useful to be able to apply the rule ``all at once''.  For instance,
if $n$ is the rank of $\G$, than a codimension-$n$ or - $(n-1)$ rule guarantees that the multiplicity of the 
corresponding component is one.

\np
Suppose that $\G$ has type $\A_4$.
In order to avoid calculating in the cohomology ring of a two-step Grassmanian when working out the 
codimension-two reduction rule, 
we use a method explained in \S\ref{sec:oneruletorulethemall} below.  
Start with $w_1=s_3s_4s_2$, $w_2=s_4s_2s_3$,
and $w=s_2s_3s_4s_2s_3s_2$, which 
have the property that $\invset_{w}=\invset_{w_1}\sqcup \invset_{w_2}$. Let $\I=\{\alpha_1,\alpha_2\}$. 
The elements $\widetilde{w}_1=s_3s_4$, $\widetilde{w}_2=s_4s_2s_3$, and $\widetilde{w}=s_2s_3s_4s_2s_3$ 
are the minimal representatives of
$w_1$, $w_2$, and $w$ in the corresponding cosets of $\W_{\PI}$, and therefore, as explained in
\S\ref{sec:oneruletorulethemall}, satisfy \eqref{eqn:Iconditions} with respect to $\I$.
For dominant weights $\mu_1=(a_1,\ldots, a_4)$, $\mu_2=(b_1,\ldots, b_4)$, and $\mu=(c_1,\ldots, c_4)$ we have

\begin{eqnarray*}
\widetilde{w}_1^{-1}\mu_1 & = & (a_1, a_2+a_3, a_4, -a_3-a_4) \\
\widetilde{w}_2^{-1}\mu_2 & = & (b_1+b_2, b_3+b_4, -b_2-b_3-b_4,b_2+b_3) \\
\widetilde{w}^{-1}\mu & = & (c_1+c_2+c_3,c_4,-c_3-c_4,-c_2).\\
\end{eqnarray*}

\np
The group $\Gtil$ is of type $\A_2$, and restriction to $\Gtil$ ignores the last two coordinates in the expressions
above, so 

\begin{equation}\label{eqn:two-step-reductions}
\left\{
\mbox{
\begin{minipage}{0.28\textwidth}
\begin{eqnarray*}
\mut_1 & = & (a_1, a_2+a_3) \\
\mut_2 & = & (b_1+b_2, b_3+b_4) \\
\mut & = & (c_1+c_2+c_3,c_4). \\
\end{eqnarray*}
\end{minipage}
}
\right.
\end{equation}

\np
The condition that $\overline{w}_1^{-1}\mu_1+\overline{w}_2^{-1}\mu_2-\overline{w}^{-1}\mu\in\Span_{\QQ}\I$ is given
by the two linear conditions

\begin{equation}\label{eqn:two-step-conditions}
\left\{
\mbox{
\begin{minipage}{0.8\textwidth}
\begin{eqnarray*}
2c_1-c_2-4c_3-2c_4 & = & (2a_1+4a_2+a_3+3a_4)+(2b_1-b_2+b_3-2b_4),\,\,\mbox{and} \\
c_1-3c_2-2c_3-c_4 & = & (a_1+2a_2-2a_3-a_4) + (b_1+2b_2+3b_3-b_4). \\
\end{eqnarray*}
\end{minipage}
}
\right.
\end{equation}

\np
This gives:

\np
{\bf Reduction rule:}
If \eqref{eqn:two-step-conditions} holds, then
$\mult_{\A_4}(\V_{\mu},\V_{\mu_1}\otimes\V_{\mu_2}) = \mult_{\A_2}(\V_{\mut},\V_{\mut_1}\otimes\V_{\mut_2})$, where
$\mut_1$, $\mut_2$, and $\mut$ are given by \eqref{eqn:two-step-reductions}.

\np
{\bf Example:} $\mu_1=(12,2,7,4)$, $\mu_2=(3,6,4,15)$, $\mu=(22,1,1,7)$, $\mut_1=(12,9)$, $\mut_2=(9,19)$, 
$\mut=(24,7)$; the multiplicity is 2.


\bpoint{A $\D_{n}$ to $\D_{n-1}$ reduction rule} 
Let $\G$ be of type $\D_{n}$ 
and let $\I=\{\alpha_2,\ldots, \alpha_n\}$.
The quotient $\Q_{n}:=\G/\PI$ is a smooth quadric hypersurface in 
$\PP^{2n-1}$.  The cohomology ring of $\Q_n$ is generated by $h$ (the class of a hyperplane section) and two
classes $a$ and $b$ of complex codimension $(n-1)$ (i.e., in the middle cohomology of $\Q_n$) satisfying the relations

\ifthenelse{\boolean{showpicture}}{%
\begin{centering}
\begin{tabular}{c}
\begin{pspicture}(0,-1.2)(4,1.2)
\pscircle[fillstyle=solid,fillcolor=gray](0,0){0.1}
\psline(0,0)(2.5,0)
\multido{\n=0.0+0.75}{4}{%
\pscircle[fillstyle=solid,fillcolor=gray](\n,0){0.1}
}
\psline(3.5,0)(4.5,0)
\SpecialCoor
\psline(!4.5 60 cos 0.75 mul add -60 sin 0.75 mul)(4.5,0)(!4.5 60 cos 0.75 mul add  60 sin 0.75 mul)
\multido{\n=3.75+0.75}{2}{%
\pscircle[fillstyle=solid,fillcolor=gray](\n,0){0.1}
}
\rput(3,0){\small \ldots}
\pscircle[fillstyle=solid,fillcolor=gray](!4.5 60 cos 0.75 mul add 60 sin 0.75 mul){0.1}
\pscircle[fillstyle=solid,fillcolor=gray](!4.5 60 cos 0.75 mul add -60 sin 0.75 mul){0.1}
\rput(0,-0.25){\tiny \oldstylenums{1}}
\rput(0.75,-0.25){\tiny \oldstylenums{2}}
\rput(1.5,-0.25){\tiny \oldstylenums{3}}
\rput(2.25,-0.25){\tiny \oldstylenums{4}}
\rput(3.75,-0.25){\tiny $n-\oldstylenums{3}$}
\rput(5.1,0){\tiny $n-\oldstylenums{2}$}
\rput(!4.5 60 cos 0.75 mul add 60 sin 0.75 mul 0.25 add){\tiny $n-\oldstylenums{1}$}
\rput(!4.5 60 cos 0.75 mul add -60 sin 0.75 mul 0.25 sub){\tiny $n$}

\end{pspicture}
\end{tabular}\\
\end{centering}
}{}

\begin{equation}\label{eqn:Qnrelations}
\rule{1cm}{0cm}h^{n-1}=a+b,ha=hb, h^na=0, a^2=b^2=\textstyle\frac{1}{2}(1-(-1)^{n})[pt], ab=\frac{1}{2}(1+(-1)^{n})[pt],
\end{equation}

\np
where $[pt]$ indicates the class of a point.  The cohomology ring of $\Q_n$ therefore has the presentation

$$
\H^{*}(\Q_{n},\ZZ) = \frac{\ZZ[h,a,b]}{(\mbox{relations in \eqref{eqn:Qnrelations}})}.
$$

\np
The integral basis for $\H^{*}(\Q_n,\ZZ)$ given by $\{h^{k}\}_{0\leq k\leq n-2}$ in codimension $\leq n-2$, 
$a$ and $b$ in codimension $n-1$, and $\{h^{k}a\}_{1\leq k\leq n-1}$ in codimensions $n$ to $2(n-1)$ is a basis of
Schubert classes in $\H^{*}(\Q_n,\ZZ)$.  We will only work out the most elementary example of a $\D_{n}$ to 
$\D_{n-1}$ reduction rule.  If $k\leq n-2$ then $h^{k}$ is the class of a Schubert cycle in $\H^{*}(\Q_{n},\ZZ)$
and the pullback to $\X$ is the class $[\Omega_{s_ks_{k-1}\cdots s_1}]$, as in the $\A_n$ case.  For $k\leq n-3$
the action of $s_k\cdots s_1$ on dominant weights is also given by the same formula as in the $\A_n$ case. 

\np
For $0\leq i,j,k\leq n-3$ with $k=i+j$, set $w_1=s_is_{i-1}\cdots s_1$, $w_2=s_js_{j-1}\cdots s_1$, and
$w=s_ks_{k-1}\cdots s_1$.  A short computation (which we omit) shows that 
$w_1^{-1}\cdot 0+w_2^{-1}\cdot 0-w^{-1}\cdot 0\in \Span_{\ZZ_{\geq 0}}\I$, and 
so $w_1$, $w_2$, and $w$ satisfy \eqref{eqn:Iconditions} with respect to $\I$.

\np
For dominant weights $\mu_1=(a_1,\cdots, a_n)$, $\mu_2=(b_1,\ldots, b_n)$ and $\mu=(c_1,\cdots, c_n)$ the
condition that $w_1^{-1}\mu_1+w_2^{-1}\mu_2-w^{-1}\mu\in\Span_{\QQ}\I$ is 

\begin{equation}\label{eqn:Dnreduction}
\rule{1.00cm}{0cm}
2\left(\sum_{r=k+1}^{n-2} c_r \right) +c_{n-1}+c_{n} = 
2\left(\sum_{r=i+1}^{n-2} a_r\right) +a_{n-1}+a_{n} + 
2\left(\sum_{r=j+1}^{n-2} b_r\right) +b_{n-1}+b_{n}. 
\end{equation}

\np
{\bf Reduction rule:} 
For any $0\leq i,j,k\leq n-3$ with $k=i+j$, if $\mu_1$, $\mu_2$, and $\mu$ satisfy 
\eqref{eqn:Dnreduction} then 
$\mult_{\D_{n}}(\V_{\smu},\V_{\smu_1}\otimes\V_{\smu_2})=\mult_{\D_{n-1}}(\V_{\smut},\V_{\smut_1}\otimes
\V_{\smut_2})$ where $\mut_1$, $\mut_2$, and $\mut$ are given by \eqref{eqn:GHreductions}.

\np
{\bf Example:} $n=5$, $i=j=1$, $k=2$, $\mu_1=(7,1,6,5,7)$, $\mu_2=(4,1,4,3,4)$, $\mu=(1,1,16,4,7)$,
$\mut_1=(8,6,5,7)$, $\mut_2=(5,4,3,4)$, $\mut=(1,17,4,7)$; the multiplicity is $514$.

\np
In order to get a $\D_{n}$ to $\D_{n-1}$ rule where the reduction formulas are different from the $\A_{n}$ case,
one only has to use deeper cohomology classes (e.g., multiplications involving $a$ or $b$).
Similar ``$\A_n$-like'' formulas hold for $\C_{n}$ to $\C_{n-1}$ and $\B_{n}$ to $\B_{n-1}$ reductions if one uses
low-codimension multiplications in $\G/\PI$ ($\I=\{\alpha_2,\ldots, \alpha_n\}$ as above), 
although the condition to check in order to apply the rule is different 
(e.g., compare \eqref{eqn:Dnreduction} and \eqref{eqn:GHconditions}).

\bpoint{A $\C_{n}$ to $\A_{n-1}$ reduction} Let $\G$ be of type $\C_n$ and 
$\I=\{\alpha_1,\ldots, \alpha_{n-1}\}$. The quotient 
$\LGn:=\G/\PI$ is the {\em Lagrangian Grassmanian}, the Grassmanian of Lagrangian $n$-planes in a $2n$-dimensional
complex vector space with a non-degenerate skew-symmetric form.

\np
Similar to the ordinary Grassmanians, the Schubert basis for $\H^{*}(\LGn,\ZZ)$ is given by classes
$\sigma_{a_1,a_2,\ldots, a_m}$ so that the corresponding partition $(a_1,a_2,\ldots,a_m)$ 
fits into an $n\times n$ box, but with the additional restriction that the partition be {\em strict}, i.e., 
that $a_1>a_2>\ldots > a_m\geq 1$
(see \cite[p.\ 29]{FP}).  
For any $a\geq 1$ set $u_{a}=s_{n+1-a}s_{n+2-a}\cdots s_{n-1}s_n$; then for any strict partition 
$n\geq a_1> a_2 > \cdots > a_m \geq 1$ the pullback of the class $\sigma_{a_1,a_2,\ldots, a_m}$ to $\H^{*}(\X,\ZZ)$ 
is the class $[\Omega_{w}]$ with $w=u_{a_m}u_{a_{m-1}}\cdots u_{a_2}u_{a_1}$.

\np
We will give only the simplest reduction rule, corresponding to the multiplication 
$\sigma_1\cdot\sigma_2=2\sigma_{3}+\sigma_{2,1}$ in cohomology.  
We must choose $w$ to be $w=u_1u_2$ (i.e., so that $[\Sigma_{w}]$ is the pullback of $\sigma_{2,1}$) 
in order to satisfy
\eqref{eqn:altIconditions}({\em ii}).
Setting $w_1=s_n$, $w_2=s_{n-1}s_n$, and $w=s_{n}s_{n-1}s_n$, then $w_1$, $w_2$, and $w$ 
satisfy \eqref{eqn:Iconditions} with respect to $\I$ (condition \eqref{eqn:Iconditions}({\em iii}) holds 
since $w_1\cdot 0 + w_2\cdot 0 - w\cdot 0 = 2\alpha_{n-1} \in \Span_{\ZZ_{\geq 0}}\I$).  The group we are reducing to
is of type $\A_{n-1}$, obtained by removing the last vertex of the Dynkin diagram for $\C_n$.

\np
If $\mu_1=(a_1,\cdots, a_n)$, $\mu_2=(b_1,\ldots, b_n)$ and $\mu=(c_1,\cdots, c_n)$ are dominant weights then 

\begin{eqnarray*}
w_1^{-1}\mu_1 & = & (a_1,a_2,\ldots,a_{n-3}, a_{n-2},a_{n-1}+2a_{n},-a_{n}) \\
w_2^{-1}\mu_2 & = & (b_1,b_2,\ldots, b_{n-3},b_{n-2}+b_{n-1},b_{n-1}+2b_{n},-b_{n-1}-b_{n}) \\
w^{-1}\mu & = & (c_1,c_2,\ldots, c_{n-3},c_{n-2}+c_{n-1}+2c_{n},c_{n-1},-c_{n-1}-c_{n}). \\
\end{eqnarray*}

\np
Restriction to $\Gtil$ ignores the last entry, so 

\begin{equation}\label{eqn:Cnreduction}
\left\{
{
\begin{array}{ccl}
\mut_1 & = & (a_1,a_2,\ldots,a_{n-3}, a_{n-2},a_{n-1}+2a_{n}) \\
\mut_2 & = & (b_1,b_2,\ldots, b_{n-3},b_{n-2}+b_{n-1},b_{n-1}+2b_{n}) \\
\mut & = & (c_1,c_2,\ldots, c_{n-3},c_{n-2}+c_{n-1}+2c_{n},c_{n-1}). \\
\end{array}
}
\right.
\end{equation}

\np
The condition that $w_1^{-1}\mu_1+w_2^{-1}\mu_2-w^{-1}\mu$ lie in $\Span_{\QQ}\I$ is 

\begin{equation}\label{eqn:Cncondition}
\sum_{r=1}^{n} r c_r - 2c_{n-1}-4c_n = 
\sum_{r=1}^{n} r a_r - 2a_n +
\sum_{r=1}^{n} r b_r - 2b_{n-1}-2b_n.
\end{equation}

\np
{\bf Reduction rule:}  If \eqref{eqn:Cncondition} holds then 
$\mult_{\C_{n}}(\V_{\smu},\V_{\smu_1}\otimes\V_{\smu_2})=\mult_{\A_{n-1}}(\V_{\smut},\V_{\smut_1}\otimes
\V_{\smut_2})$ 
where $\mut_1$, $\mut_2$, and $\mut$ are given by \eqref{eqn:Cnreduction}.



\np
{\bf Example:}  $n=5$, $\mu_1=(8,4,3,1,3)$, $\mu_2=(3,2,1,6,1)$, $\mu=(6,6,14,1,1)$, $\mut_1=(8,4,3,7)$, 
$\mut_2=(3,2,7,8)$, and $\mut=(6,6,17,1)$; the multiplicity is $31$.



\np
{\bf Remark on saturation.} 
If $(\mu_1,\ldots,\mu_k,\mu)$ is an integral point of $\LRC(k)$ it does not necessarily imply that
$\V_{\mu}$ is a component of $\V_{\mu_1}\otimes\cdots\otimes\V_{\mu_k}$. 
The problem of determining 
the integral points for which this implication does hold is known as the saturation problem.  
For any integral point of $\LRC(k)$ it is known that the implication holds for some positive multiple of that point, 
and that the multiple can be bounded by a constant depending only on $\G$. 
The cone $\LRC(k)$ (respectively a face $\F$ of $\LRC(k)$) is called {\em saturated} 
if the implication holds for every integral point in the cone (respectively on the face).
In type $\A$, all cones are saturated by the theorem of Knutson-Tao \cite[p.\ 1084]{KT}.  If $\F$ is a regular face
such that the corresponding reduction rule reduces to a group of type $\A$, as in the example above, then
the reduction theorem and the result of Knutson-Tao imply that $\F$ is saturated.

\bpoint{A rule for producing reduction rules}  \label{sec:oneruletorulethemall}
Suppose that $w_1$,\ldots, $w_k$, and $w$ are elements of $\W$
such that 

\begin{equation}\label{eqn:liningup}
\invset_{w} = \bigsqcup_{i=1}^{k} \invset_{w_i},
\end{equation}

\np
i.e., $\invset_{w}$ is the disjoint union of $\invset_{w_1}$ through $\invset_{w_k}$.  In the classical cases
one can check that \eqref{eqn:liningup} implies that $\cap_{i=1}^{k}[\Omega_{w_i}]\cdot[\X_{w}]=1$, 
and an argument proving this for general semisimple $\G$ will appear in \cite{kp}.  
If $\I'$ is the empty set (so $\P_{\I'}=\B$ and $\W_{\P_{\I'}}=\{e\}$) then $w_1$,\ldots, $w_k$, and $w$ 
satisfy \eqref{eqn:Iconditions} with respect to $\I'$ (condition \eqref{eqn:Iconditions}({\em iii}) follows from
\eqref{eqn:liningup} and \eqref{eqn:rootsum}).
Thus $w_1$,\ldots, $w_k$, and $w$ describe a codimension-$n$ regular face of $\LRC(k)$ 
and a corresponding codimension-$n$ reduction rule, where $n$ is the rank of $\G$.

\np
Furthermore, for any subset $\I$ of the simple roots, if we set $\widetilde{w}_1$,\ldots, $\widetilde{w}_k$, and 
$\widetilde{w}$ to be the shortest elements in the cosets $w_1\W_{\PI}$,\ldots, $w_k\W_{\PI}$, and $w\W_{\PI}$ 
respectively, then \cite[Lemma 7.1.3]{dr1} shows that $\widetilde{w}_1$,\ldots, $\widetilde{w}_k$ and $\widetilde{w}$ 
satisfy
\eqref{eqn:Iconditions} with respect to $\I$, yielding a codimension $n-|\I|$ face of $\LRC(k)$ and a corresponding
reduction rule.  I.e., the elements $w_1$,\ldots, $w_k$, and $w$ give a family of reduction rules, one for each
subset $\I$ of simple roots.  I do not know if all regular faces arise via this procedure.  

\np
Any face containing a regular face is itself regular, and of course, the codimension $n-|\I|$ faces above are simply
all the faces containing the codimension-$n$ face corresponding to $w_1$,\ldots, $w_k$, and $w$.   The question as to
whether all regular faces arise via the procedure above is therefore equivalent to the question as to whether
every regular face contains a regular codimension-$n$ face.

\section{Further remarks}

\bpoint{GIT interpretation of the reduction theorem}  Suppose that $\F$ is a regular face of $\LRC(k)$, and 
let $\I$, $w_1$, \ldots, $w_k$, and $w$ be the data parametrizing $\F$ given by Theorem \ref{thm:Ressayre}. 
Let $\psi\colon \Xtil^{k+1}\longrightarrow\X^{k+1}$ be the embedding given in the reduction theorem. 
For any point $(\mu_1,\ldots, \mu_k,\mu)$ of $\F$ in the strictly dominant locus, the line bundle
$\L:=\L_{-\wo\mu_1}\bt\cdots\bt\L_{-\wo\mu_k}$ $\bt\L_{\mu}$ 
is very ample on $\X^{k+1}$, and hence
its pullback $\Ltil:=\psi^{*}\L$ is very ample on $\Xtil^{k+1}$.
For all $m\geq 0$ the reduction theorem implies that pullback by $\psi$ induces an isomorphism
$\psi^{*}\colon\H^{0}(\Xtil^{k+1},\L^{m})^{\G}\stackrel{\sim}{\longrightarrow}\H^{0}(\Xtil^{k+1},\Ltil^{m})^{\Gtil}$. 

\np
The $\Gtil$-equivariant embedding $\psi$ induces a map of GIT quotients 
$\displaystyle\Xtil^{k+1}\quot_{{\Ltil}}\Gtil\longrightarrow \X^{k+1}\quot_{\L}\G$, and the equality of pullbacks
above for all $m\geq 0$ implies that this map is an isomorphism.

\np
\bpoint{Reduction to $\LRC_{\Gtil}(k)$}
If $\F$ is a regular face of $\LRC(k)$, 
and $\Gtil$ the corresponding group provided by the theorem, then reduction gives a map from $\F$ to $\LRC_{\Gtil}(k)$, 
the Littlewood-Richardson cone of $\Gtil$.
In this section we prove some basic results about this reduction map.

\np
Recall that for any polyhedral cone $\mathcal{C}$ 
in a vector space $\E$ 
every point $p\in\mathcal{C}$ lies on the relative interior of a unique face.  The dimension of this face
is the same as the dimension of the subspace of $\E$ spanned by the set
$\left\{\varepsilon\in\E \st p\pm \varepsilon\in\mathcal{C}\right\}.$

\tpoint{Proposition}\label{prop:reduction}
Suppose that $\F$ is a regular face of codimension $r$, 
that $p=(\mu_1,\ldots, \mu_k$, $\mu)$  is a point of $\F$ in the strictly dominant locus, and that $p$ lies on the 
relative interior of a face of $\LRC(k)$ of codimension $r'$.
Then the image
of $p$ under the reduction map $\F\longrightarrow\LRC_{\Gtil}(k)$ lies on the relative interior of a regular face of codimension
$r'-r$.

\np
\bpf
Let $(\mut_1,\mup_1)$, \ldots, $(\mut_k,\mup_k)$, and $(\mut,\mup)$ be the restrictions of $w_1^{-1}\mu_1$,\ldots
$w_k^{-1}\mu_k$, and $w^{-1}\mu$ respectively to $\gttil$ and $\ga$ under the splitting $\gt=\gttil\oplus\ga$ from \S\ref{sec:GtilandGI},
so that 
$\ptil:=(\mut_1,\ldots,\mut_k,\mut)$ is the image of $p$ under the reduction map.   
By the discussion at the end of \S\ref{sec:GtilandGI}, $\pu$ is strictly dominant, and so we only need to check
the statement on codimension.  
Write

$$p=
\left({(\mut_1,\mup_1),\ldots, (\mut_k,\mup_k),(\mut,\mup)\rule{0cm}{0.4cm}}\right),$$

\np
meaning that we have acted by $w_i^{-1}$ (or $w^{-1}$) and applied the splitting to each entry.
Let 
$\ept:=(\ept_1,\ldots, \ept_k,\ept_{k+1})$ be a tuple with each $\ept_i\in\gttil^{*}$, 
$\epp=(\epp_1,\ldots,\epp_k,\epp_{k+1})$ be a tuple with each $\epp_i\in \ga^{*}$, and

$$p\pm (\ept,\epp):= \left({(\mut_1\pm\ept_1,\mup_1\pm\epp_1),\ldots, 
(\mut_k\pm\ept_k,\mup_k\pm\epp_k),
(\mut\pm\ept_{k+1},\mup\pm\epp_{k+1})\rule{0cm}{0.4cm}}\right). $$

\np
The vector space map underlying the reduction map sends $p\pm(\ept,\epp)$ to $\ptil\pm \ept$. We want to study
which $\ept$ such that $\ptil\pm\ept\in \LRC_{\Gtil}(k)$ can be realized as the image of $(\ept,\epp)$ such that
$p\pm(\ept,\epp)\in\LRC(k)$.  We can make the following simplifying assumptions: ({\em i}) since both $\LRC(k)$ and
$\LRC_{\Gtil}(k)$ are rational cones, we may restrict to rational $\ept$ and $\epp$. ({\em ii}) since it is only the
dimension of the vector space spanned by $\ept$ (respectively $(\ept,\epp)$) that matters, we may scale these vectors
and assume they are arbitrarily small.  In particular, since $p$ is strictly dominant, we may (after scaling 
$(\ept,\epp)$), assume that both $p\pm(\ept,\epp)$ are dominant.

\np
Since $p$ lies on the face $\F$, in order for $p\pm(\ept,\epp)$ to be in $\LRC(k)$ a necessary condition is that
$p\pm(\ept,\epp)$ satisfy the linear conditions defining $\F$.  In these coordinates, the condition is simply that
$\sum \epp_i=0\in\ga^{*}$.  If this condition holds, (and since $p\pm(\ept,\epp)$ are dominant) we may apply the
reduction rule.  Scaling by some positive integer $m$ so that all weights are integral, the reduction rule says that

$$
\mult_{\G}(
\V_{m(\mut\pm\ept_{k+1},\mup\pm\epp_{k+1})},
\V_{m(\mut_1\pm\ept_1,\mup_1\pm\epp_1)}
\otimes\cdots\otimes
\V_{m(\mut_k\pm\ept_k,\mup_k\pm\epp_k)})
=\rule{5cm}{0cm}$$
$$
\rule{7.5cm}{0cm}
\mult_{\Gtil}(
\V_{m(\mut\pm\ept_{k+1})},
\V_{m(\mut_1\pm\ept_1)}
\otimes\cdots\otimes
\V_{m(\mut_k\pm\ept_k,)}).
$$

\np
Thus (subject to the simplifying assumptions above),  $p\pm(\ept,\epp)\in\LRC(k)$ if and only if $\sum \epp_i=0\in\ga^{*}$
and $\ptil\pm\ept\in\LRC_{\Gtil}(k)$. In particular this shows that (up to scaling) all $\ept$ such that
$\ptil\pm\ept\in\LRC_{\Gtil}(k)$ may be realized, and that the kernel of the map $(\ept,\epp)\longrightarrow\ept$
has dimension $k\dim_{\CC}(\ga^{*})=kr$.  Counting dimensions then gives the proposition. \epf

\np
Here are some immediate corollaries.
First, 
the proposition implies the result mentioned in \S\ref{sec:codim-two}.

\tpoint{Corollary}\label{cor:codimr-restriction} The reduction rule corresponding to a regular face of codimension
$r$ can be obtained as a succession of $r$ codimension-one reduction rules.

\np
\bpf Suppose that $\F$ is a regular face of codimension $r$, then $\F$ is contained in a codimension $1$ face
$\F'$ which must also be regular.   Let $\Gtil'$ be the group corresponding to $\F'$, then by Proposition
\ref{prop:reduction}  the image of $\F$ under the codimension-one reduction map 
$\F'\longrightarrow\LRC(k)_{\Gtil'}$ is a regular face of codimension $r-1$.  Continuing inductively 
we obtain a succession of $r$ codimension-one reduction rules.  
What remains is to check that the composition of these rules is the same rule as the codimension-$r$ rule obtained
from the face $\F$.  We briefly sketch how to produce at least one factorization such that this holds.

\np
Suppose that the face $\F$ is determined by the data
$\I$, $w_1$,\ldots, $w_k$, and $w$ as in Theorem \ref{thm:Ressayre}.  Let $\alpha_j\in\I$ be any simple root, 
and $\P_{j}$ the parabolic obtained by inverting $\alpha_j$.  Let $\widetilde{w}_1$,\ldots, $\widetilde{w}_k$, and
$\widetilde{w}$ be the minimal representatives in the cosets
$w_1\W_{\P_{j}}$, \ldots
$w_k\W_{\P_{j}}$, and
$w\W_{\P_{j}}$ respectively, and let $u_1$,\ldots, $u_k$, and $u\in \W_{\P_{j}}$ be the unique
elements such that $w_1=\widetilde{w}_1u_1$,\ldots, $w_k=\widetilde{w}_ku_k$, and $w=\widetilde{w}u$.  Then
similarly to the proof of \cite[Lemma 7.1.3]{dr1} one can check that $\widetilde{w}_1$,\ldots, 
$\widetilde{w}_k$, and $\widetilde{w}$ satisfy conditions \eqref{eqn:Iconditions} with respect to $\{\alpha_j\}$ 
and so define an codimension-one face $\F'$.  Furthermore, $u_1$,\ldots, $u_k$, and $u$ satisfy
\eqref{eqn:Iconditions} with respect to $\I\setminus\{\alpha_j\}$ in the group $\Gtil'$, and parametrize the
regular face corresponding to the image of $\F$ in $\LRC(k)_{\Gtil'}$.  
The corresponding codimension-one reduction rule is computed in coordinates (as in the examples
above) by writing $\widetilde{w}_1^{-1}\mu_1$,\ldots, $\widetilde{w}_k^{-1}\mu_k$, and 
$\widetilde{w}^{-1}\mu$ in the basis of 
fundamental weights and dropping the $j$-th coordinate.  This is the same as 
writing $w_1^{-1}\mu_1$,\ldots, $w_k^{-1}\mu_k$, and $w^{-1}\mu$ in the basis of fundamental weights,
dropping the $j$-th coordinate, and then applying $u_1$,\ldots, $u_k$, and $u$ to the result.  
This shows that the composition of the codimension-one and codimension-$(r-1)$ rule is equal to the codimension-$r$
rule, and by induction that the composition of the succession of $r$ codimension-one rules is equal to the original
codimension-$r$ rule. \epf

\np
Second, by taking a point $p$ in the relative interior of $\F$ we obtain:

\tpoint{Corollary}\label{cor:restriction} 
The image of the reduction map is a full dimensional subcone of $\LRC_{\Gtil}(k)$.

\np
This reduction map is not surjective in general, and it would be interesting to know how to characterize
the image.

\bigskip
\bigskip
\np
{\small
Department of Mathematics and Statistics, \\
Queen's University, Kingston, Ontario, Canada, K7L 3N6.\\
{\em E-mail address}: {\tt mikeroth@mast.queensu.ca}}

\end{document}